\documentclass[11pt,reqno]{amsart}

\usepackage{amsmath}
 \usepackage{amssymb}
 \usepackage{color}

\newcommand{\mysection}[1]{\section{#1}
      \setcounter{equation}{0}}

\newcommand{\nlimsup}{\operatornamewithlimits{\overline{lim}}}

\newtheorem{theorem}{Theorem}[section]
\newtheorem{lemma}[theorem]{Lemma}

\newtheorem{corollary}[theorem]{Corollary} 

\theoremstyle{definition}
\newtheorem{assumption}{Assumption}[section]
\newtheorem{definition}{Definition}[section]
\theoremstyle{remark}
\newtheorem{remark}{Remark}[section]

\newcommand{\tr}{\text{\rm tr}\,}

\newcommand\bbeta{\text{\raise-.2ex\hbox{$\bm{\beta}$}}}
\newcommand\dist{{\rm dist}\,}

\makeatletter
 \def\dashint{%  
 \operatorname%
 {\,\,\text{\bf--}\kern-.98em\DOTSI\intop\ilimits@\!\!}}
 \makeatother

\newcommand\bR{\mathbb{R}}

\newcommand\bS{\mathbb{S}}

\newcommand\cR{\mathcal{R}}

\newcommand\cH{\mathcal{H}}

\newcommand\sfu{{\sf u}}

\renewcommand\({{\rm(}}

\renewcommand\){{\rm)}}

\newcommand\supinf{\operatornamewithlimits{sup\,\,\,inf}}
\newcommand\supsup{\operatornamewithlimits{sup\,\,\,sup}}

\newcommand\loc{\textnormal{loc}}
\newcommand\esssup{\operatornamewithlimits{ess\,sup\,}}

\begin{document}

\title[On Isaacs equations]
{Uniqueness for $L_{p}$-viscosity solutions for
uniformly parabolic Isaacs equations
with measurable lower order terms}

\author{N.V. Krylov}

\email{nkrylov@umn.edu}
\address{127 Vincent Hall, University of Minnesota,
 Minneapolis, MN, 55455}

\keywords{Fully nonlinear equations,
viscosity solutions, Isaacs equations}

\renewcommand{\subjclassname}{%
\textup{2010} Mathematics Subject Classification}

\subjclass[]{35K55, 35B65}

\begin{abstract}
 In this article we present several results
concerning uniqueness of $C$-viscosity
and $L_{p}$-viscosity solutions 
for fully nonlinear parabolic equations.
In case of the Isaacs equations we allow
lower order terms to have just measurable bounded 
coefficients. Higher-order coefficients
are assumed to be H\"older continuous
in $x$ with exponent slightly less than $1/2$.
This case is treated by using stability
of maximal and minimal $L_{p}$-viscosity solutions.
 
\end{abstract}

\maketitle

\mysection{Introduction}
                                            \label{section 10.25.1}

For a real-valued measurable function   $H(\sfu,t,x)$,    \smallskip
$$
\sfu=(\sfu',\sfu''),\quad
\sfu'=\big(\sfu'_{0},\sfu'_{1},...,\sfu'_{d}\big) \in\bR^{d+1},\quad 
\sfu''\in\bS,\quad
(t,x)\in\bR^{d+1}, 
\medskip$$ 
where $\bS$ is the
set of symmetric $d\times d$ matrices, and
sufficiently regular functions $v(t,x)$ we set  \smallskip
$$
H[v](t,x)=H\big(v(t, x),D v(t,x),D^{2}v(t, x),t, x\big),
\medskip$$
and
we will be dealing with 
the parabolic  equations   \smallskip
\begin{equation}
                                                \label{7.29.10}
\partial_t v(t,x)+H[v](t,x)=0
\smallskip\end{equation}

\noindent
in subsets of $[0,T)\times \bR^d $, where $T\in(0,\infty)$
is fixed. Above  \smallskip
$$
\bR^{d}=\big\{x=(x^{1},...,x^{d}):x^{1},...,x^{d}\in \bR\big\},
\medskip$$
$$
\partial_t=\frac{\partial}{\partial t},\quad
 D^{2}u=(D_{ij}u),\quad Du=(D_{i}u),\quad
D_{i}=\frac{\partial}{\partial x^{i}},
\quad D_{ij}=D_{i}D_{j}.
\medskip$$

If $  R\in(0,\infty)$ and $(t,x)\in\bR^{d+1}$, then   \smallskip
$$
B_{R}=\big\{x\in\bR^{d}:|x|<R\big\},\quad B_{R}(x)=x+B_{R},
\medskip$$
$$
C_{ R}=[0,R^{2})\times B_{R},\quad
C_{R}(t,x)=(t,x)+C_{R}.
\medskip$$

We also take a bounded domain $\Omega\subset\bR^{d}$
of class $C^{1,1}$ and set   \smallskip    
$$
\Pi=[0,T)\times\Omega,\quad\partial'\Pi=\bar\Pi 
\setminus\big(\{0\}\times\bar\Omega\big)
\medskip$$
 
 \begin{remark}
We assumed that $\Omega\in C^{1,1}$ just to be able to refer
to the results available at this moment, but actually
much less is needed for our Theorems
\ref{theorem 3,13,1},  \ref{theorem 10.9.1},  \ref{theorem 1,31,3},
\ref{theorem 3,11,1},  \ref{theorem 8,12,1}, and \ref{theorem 2,3,1} 
to hold. For instance the exterior cone condition would
suffice.

\end{remark}
We will be dealing with viscosity solutions of 
\eqref{7.29.10} in $\Pi$. The following definition
is  taken from \cite{CKS00}
and has the same spirit as in  \cite{CIL}.

\begin{definition}
                                             \label{definition 1,29.1}  
  For each
choice of ``regularity'' class $\cR=C $ or $\cR=L_{p}$
we say that $u$ is an $\cR$-viscosity subsolution
\index{$\cR$-viscosity solution}%
 of 
\eqref{7.29.10} in $\Pi$ provided that $u$ is continuous
in $\Pi$ and, for any $\bar C_{r}(t_{0},x_{0})\subset\Pi$  and any function
$\phi$, that is   continuous in $C_{r}(t_{0},x_{0})$ and
whose generalized derivatives satisfy
$\partial_{t}\phi,D\phi,D^{2}\phi\in\cR\big(C_{r}(t_{0},x_{0})\big)$, and is such that $u-\phi$ attains its
maximum over $C_{r}(t_{0},x_{0})$ at $(t_{0},x_{0})$, we have \medskip
\begin{equation}
                                               \label{1,29,2}
 \lim _{ \rho\downarrow0}\esssup_{C_{\rho}(t_{0},x_{0}) }
\big[\partial_{t}\phi(t,x) +
H\big(u(t,x),D\phi(t,x),D^{2}\phi(t,x),t,x\big)\big]\geq 0.
\medskip\end{equation}

In a natural way one defines $\cR$-viscosity supersolutions
and calls a function an $\cR$-viscosity solution if it is an
$\cR$-viscosity supersolution and an $\cR$-viscosity subsolution. 
\end{definition}

Note that $C_{r}(t_{0},x_{0})$ contains $\big\{(t,x):
t=t_{0}, |x-x_{0}|<r\big\}$, which is part of its boundary. Therefore,
the conditions like $D^{2}\phi\in C\big(C_{r}(t_{0},x_{0})\big)$ mean that
the second-order derivatives of $\phi$ are continuous
up to this part of the boundary. 

In Section \ref{section 3,13,1} we discuss uniqueness
of $C$-viscosity solutions for general equations when
$H$ is Lipschitz continuous with respect to $\sfu$.
The result we obtain is crucial for proving
uniqueness of $L_{p}$-viscosity solutions in 
Section \ref{section 3,7,1}
for the Isaacs equations with measurable lower order terms.
The proof of the main result in Section \ref{section 3,13,1}
 hinges on
Lemma \ref{lemma 2,4,2}, whose rather long proof is given in
  Section \ref{section 02,1,1}. Section \ref{section 12,2.1}
concentrates on the extremal $L_{p}$-viscosity
solutions, their existence and stability.
Precisely the stability of $L_{p}$-viscosity
minimal and maximal
solutions is used in 
Section \ref{section 3,7,1}.

\mysection{Uniqueness of $C$-viscosity
solutions of parabolic equations}
                                   \label{section 3,13,1}
Fix some constants $\delta\in(0,1]$, $K_{0} \in(0,\infty)$,
and set
$$
\bS_{\delta}=\{a\in\bS: \delta^{-1}|\lambda|^{2}\geq a^{ij}\lambda^{i}\lambda^{j}\geq \delta |\lambda|^{2}\quad\forall \lambda\in\bR^{d}\}.
$$

In the following assumption there is 
$\gamma=\gamma(d,\delta)\in(1/4,1/2)$ which we specify
in the following way. In Lemma 6.3 of \cite{Kr_17_1}
$$
\kappa(d,\delta,p)\in(1,2)
$$
(close to $1$) is defined.  We take $\bar \delta=\bar \delta(d,\delta)
\in(0,\delta)$
(close to $0$) from
Theorem 4.1 of \cite{Kr_17_1} and set  \vspace{5pt}
$$
\kappa=\kappa(d,\bar\delta,d+3),\quad
 \gamma=\frac{7-3\kappa}{12-4\kappa}\quad(\in(1/4,1/2).
$$

\begin{assumption}
                                   \label{assumption 3,13,1}
(i) The function $H(\sfu,t,x)$ is a {\em continuous\/} function
of $(\sfu,t,x)$ and is Lipschitz continuous with respect to $\sfu$
with Lipschitz constant $K_{0}$.

(ii) At all points of differentiability of $H$ with respect to $\sfu''$
we have $D_{\sfu''}H\in\bS_{\delta}$. 

(iii) For   all values of the arguments we have   \vspace{5pt}
\begin{equation}
                                                   \label{3,15,1}
\big|H(\sfu,t,x)-H(\sfu,t,y)\big|\leq K_{0}|x-y|^{\gamma}|\sfu''|
+\big(|\sfu'|+1\big)\omega\big(|x-y|\big),
\end{equation}  \vspace{5pt}
where $\omega(\tau)$, $\tau\geq0$, is a continuous functions
vanishing at the origin.

(iv) $ \sup \big\{|H(0,t,x)|:(t,x)\in\bR^{d+1}\big\}=:\bar H<\infty $.

(iv) We are given a $g\in C(\overline{\partial'\Pi})$.  

\end{assumption}

\begin{theorem}
                                           \label{theorem 3,13,1}
Under the above assumptions there exists
a unique $v\in C(\bar\Pi)$ which is a $C$-viscosity
solution of \eqref{7.29.10} in $\Pi$ with boundary
condition $v=g$ on $\partial'\Pi$.
Furthermore, there exists a constant   $N\in(0,\infty)$
such that for any $\rho>0$  
\begin{equation}
                                           \label{3,16,6}
\|v\|_{C^{\kappa}(\Pi^{\rho})}\leq N\rho^{-\kappa},
\vspace{5pt}\end{equation}
where $\Pi^{\rho}=[0,T-\rho^{2})\times\Omega^{\rho}$,
$\Omega^{\rho}=\{x :\rho_{\Omega}(x)
>\rho\}$, $\rho_{\Omega}(x)=\dist(x, \Omega^{c})$.

\end{theorem}

\begin{remark}
                                             \label{remark 3,23,1}
The assumptions of Theorem \ref{theorem 3,13,1}
are almost identical to the assumptions made in the elliptic case
in \cite{Tr_88}, that, to the best of our knowledge, 
provides  the most general result  to date  concerning
the uniqueness of  $C$-viscosity solutions for the 
uniformly elliptic case
(see Remark 3.1  there).    Our Theorem \ref{theorem 3,13,1}
is a parabolic counterpart of Trudinger's result 
from \cite{Tr_88}.

In the parabolic case the
uniqueness of $L_{p}$-viscosity solutions
is proved in Lemma 6.2 of \cite{CKS00}
when $H$ is independent of $(t,x)$.
In the case of the Isaacs equations, under the assumptions
on the coefficients guaranteeing that our assumptions are satisfied   
as well,
the statement about the uniqueness of $C$-viscosity solutions
is found in Theorem 9.3 of  \cite{CKS00}. However, this statement
is not provided with a proof with the excuse that its proof is similar 
to the one known in the elliptic case.

One of the features of our proof is that it also
allows one to establish an algebraic rate of convergence
of numerical approximations (see \cite{Kr15.1}).

\end{remark}

The proof of Theorem \ref{theorem 3,13,1} is based on a few
auxiliary results. Denote
$$
[\sfu']=(\sfu'_{1},...,\sfu'_{d}).
$$

\begin{remark}
                                      \label{remark 3,15,4}
It is easy to see that if $v(t,x)$ is a $C$-viscosity
subsolution of \eqref{7.29.10} in $\Pi$, then, for any 
constant $c$, the functions
$w(t,x):=e^{ct}v(t,x)$ is a  $C$-viscosity subsolution of
$$
\partial_{t}w+H^{c}[w]=0 
$$
in $\Pi$, where 
$$
H^{c}(\sfu,t,x):=e^{ct}H(e^{-ct}\sfu,t,x)-c\sfu'_{0}.
\vspace{5pt}$$
The function $H^{c}(u,t,x)$ has the same Lipschitz constant
with respect to $\big([\sfu'],\sfu''\big)$ and its derivative with
respect to $\sfu'_{0}$, wherever it exists,
is   \smallskip
$$
D_{\sfu'_{0}}H(e^{-ct}\sfu,t,x)-c\leq K_{0}-c.
\smallskip$$
If we take $c=K_{0}+1$ and redefine $H^{c}$ for $t\not\in[0,T]$
as its value at the closest end point of $[0,T]$,
then $H^{c}$ will satisfy all assumptions of Theorem 
\ref{theorem 3,13,1} with $K_{0}+1$ in place of $K_{0}$
and additionally satisfy $D_{\sfu'_{0}}H_{c}\leq-1$.

That is why without loss of generality we suppose that
not only Assumption \ref{assumption 3,13,1}
is satisfied but also for all values of arguments
\begin{equation}
                                                    \label{3,15,5}
D_{\sfu'_{0}}H \leq-1
\end{equation}
wherever the left-hand side exists.
\end{remark}

Below we suppose that the above assumptions are satisfied,   
take the convex positive homogeneous of degree one function
$P(\sfu'')$ on $\bS$ from Theorem \ref{theorem 10.9.1},
and set $P[u](t,x)=P\big(D^{2}u(t,x)\big)$.
Recall that $\kappa$ is introduced before Assumption \ref{assumption 3,13,1}.

\begin{lemma}
                                   \label{lemma 2,4,1}

\(i\,\) For any $K\geq1$ each of the equations  \smallskip
\begin{equation}
                                                 \label{3.30.6}
\partial_{t}u_{ K}+\max\big(H[u_{  K}],  P[  u
_{ K}]- K\big)=0,
\smallskip\smallskip\end{equation}
\begin{equation}
                                                 \label{3.30.06}
\partial_{t}u_{-K}+\min\big(H[u_{-K}],-P[- u
_{- K}]+K\big)=0,
\smallskip\smallskip\end{equation}
in $\Pi$ \(a.e.\) with boundary condition $u_{\pm K}=g$ on $\partial' \Pi$
has a unique solution 
 $u_{\pm K}\in W^{1,2}_{p,\loc}(\Pi^{\rho})\cap C(\bar \Pi)$
for any $p\geq1$ and $\rho>0$. 

\(ii\,\) We have $u_{-K}\leq u_{K}$ and,
as $K$  increases, $u_{-K}$ increase and $u_{K}$ decrease.

\(iii\,\) The family $\{u_{-K},u_{K}:K\geq 1\}$ is equicontinuous
in $\bar \Pi$.

\(iv\,\) There exists a constant   $N\in(0,\infty)$
such that for any $\rho>0$ and $K\geq 1$
\begin{equation}
                                           \label{3,16,7}
\|u_{K},u_{-K}\|_{C^{\kappa}(\Pi^{\rho})}\leq N\rho^{-\kappa}.
\end{equation}

\end{lemma}

The existence part in assertion (i) for the sign +  
follows from Theorem
\ref{theorem 10.9.1}
which holds under more general assumptions than the ones
imposed here.
 For the sign $-$ it suffices to replace
$H(\sfu,t,x)$ with $-H(-\sfu,t,x)$.
Uniqueness and assertion (ii)
are direct consequences of the maximum principle.

Assertion (iii) for $u_{K}$ follows from the linear theory
and the observation that   
$$
\big|\max\big(H(\sfu',0,t,x),-K\big)\big|\leq\big|H(\sfu',0,t,x)\big|\leq 
K_{0}|\sfu'|+ \bar H.
 $$
Indeed, for any $K$,
 there exist $\bS_{\bar\delta}$-valued $a$, $\bR^{d}$-valued $b$,
and real-valued $c\geq0$ and $f$ such that
$$
\partial_{t}u_{K}+a^{ij}D_{ij}u_{K}+b^{i}D_{i}u_{K}-cu_{K}+f=0
$$
and $|b|\leq K_{0}$, $c\leq K_{0}$, $|f|\leq\bar H$.
For $u_{-K}$ the argument is similar.
 
Assertion (iv) follows from Theorem 2.1 of \cite{Kr_17_1}.

\begin{lemma}
                                                 \label{lemma 2,4,2}
Under the assumptions of this section,
for $K\to\infty$ we have $|u_{K}-u_{-K}|\to0$
uniformly in $\bar\Pi$.
\end{lemma}

This lemma is proved in Section \ref{section 02,1,1}.

{\bf Proof of Theorem \ref{theorem 3,13,1}}.
First we prove uniqueness. 
Introduce $\psi\in C^{2}(\bR^{d})$ as a global barrier for $\Omega$,
that is, in $\Omega$ we have $\psi\geq1$ and   \smallskip
$$
a_{ij}D_{ij}\psi+b_{i}D_{i} \psi \leq-1
\smallskip$$
for any $(a_{ij})\in \bS_{\bar \delta}, \big|(b_{i})\big|\leq K_{0}$.
Such a $\psi$ can be found in the form
$  \cosh \mu R-\cosh\mu|  x|$ for sufficiently large $\mu$
and $R$.

Then we take and fix a radially symmetric
with respect to $x$ function $\zeta=\zeta(t,x)$ of class
$ C^{\infty}_{0}(\bR^{d+1})$ with support in $(-1,0)\times B_{1}$
and unit integral.
For $\varepsilon>0$ we define $\zeta_{\varepsilon}(t,x)
=\varepsilon^{-d-2}\zeta(\varepsilon^{-2}t,\varepsilon^{-1}x)$
and for locally summable $u(t,x)$ introduce  
\begin{equation}
                                                   \label{5.3.2}
u^{(\varepsilon)}(t,x)=u(t,x)*\zeta_{\varepsilon}(t,x).
\smallskip\end{equation}

Let $\Omega_{n}$, $n=2,3,...$, be a sequence of 
strictly expanding smooth
domains whose union is $\Omega$ and set $\Pi_{n}=[0, T(1-1/n))
\times\Omega_{n}$. Then for any $n_{0}=3,4,...$ and all sufficiently small
$\varepsilon>0$
 \smallskip
$$
\xi_{\varepsilon, K} :=\partial_{t}u_{  K}^{(\varepsilon)}
+\max\big(H\big[u_{ K}^{(\varepsilon)}\big],
 P\big[u_{  K}^{(\varepsilon)}\big]-K\big) 
\smallskip\smallskip$$
  is well defined in $\Pi_{n_{0}}$. Since the second-order derivatives  
with respect to $x$ and the first derivative 
with respect to $t$
of $u_{ K}$   are     in $L_{p}(\Pi^{\rho})$  
for any $p$ and $\rho$,
we have $ \xi_{\varepsilon, K} 
\to0$ as $\varepsilon\downarrow0$
in any $L_{p}(\Pi_{n_{0}})$ for any $K$
and any $p>1$. Furthermore,
$ \xi_{\varepsilon,K} $ are
 continuous {\em because $H(\sfu,t,x)$
is continuous\/}. Therefore, there exist smooth functions  
$ \zeta_{\varepsilon,K} $ on $\bar \Pi_{n_{0}}$ such that
$$
-\varepsilon\leq \xi_{\varepsilon,K}+
\zeta_{\varepsilon,K}\leq0 
$$
in $\Pi_{n_{0}}$ for all small $\varepsilon>0$.

Since $\Omega_{n_{0}}$ is smooth, by Theorem 1.1
of \cite{DKL_12} there exists a unique
 $
w_{\varepsilon,K}\in \bigcap_{ 
p>1}W^{1, 2}_{p }(\Pi_{n_{0}} ) 
 $
 satisfying 
$$
\partial_{t}w_{\varepsilon,K}+
\sup_{\substack{a\in\bS_{\bar\delta},|b|\leq K_{0}\\0\leq
 c \leq K_{0}}}\big[
a_{ij}D_{ij}w_{\varepsilon,K}+b_{i}D_{i}w_{\varepsilon,K}
-cw_{\varepsilon,K}\big]=
\zeta_{\varepsilon,K}
\smallskip\smallskip$$
 in $\Pi_{n_{0}} $ \(a.e.\)  and such
that $w_{\varepsilon,K}=0$ on $\partial'\Pi_{n_{0}} $. By the maximum principle
such $w_{\varepsilon,K}$ is unique. Then owing to the continuity
of $\zeta_{\varepsilon,K}$, for any $\varepsilon$ and $K$,
for all sufficiently small $\beta>0$, we have   
\begin{equation}
                                             \label{4.28.1}
\partial_{t}w^{(\beta)}_{\varepsilon,K}+
\sup_{\substack{a\in\bS_{\bar\delta},|b|\leq K_{0}\\0\leq c
\leq K_{0}}}\big[
a_{ij}D_{ij}w^{(\beta)}_{\varepsilon,K}+
b_{i}D_{i}w^{(\beta)}_{\varepsilon,K}
-cw^{(\beta)}_{\varepsilon,K}\big]
\leq
\zeta^{(\beta)}_{\varepsilon,K}\leq \zeta_{\varepsilon,K} +\varepsilon
\end{equation}
in $\Pi_{n_{0}-1}$.

 For $w^{\beta}_{\varepsilon,K}(t,x):=w^{(\beta)}_{\varepsilon,K}
+\varepsilon(T-t)$ now, obviously,   \smallskip
$$
\partial_{t}\big(u_{K}^{(\varepsilon)}+w^{\beta}_{\varepsilon,K}\big)+
\max\big(H\big[u_{K}^{(\varepsilon)}+w^{ \beta}_{\varepsilon,K}\big],
P\big[u_{K}^{(\varepsilon)}+w^{\beta}_{\varepsilon,K}\big]-K\big) 
\vspace{5pt}$$
$$
\leq \partial_{t}u_{K}^{(\varepsilon)}
+\max\big(H\big[u_{K}^{(\varepsilon)}\big],
P\big[u_{K}^{(\varepsilon)}\big]-K\big)+
\partial_{t}w^{ \beta}_{\varepsilon,K}
\vspace{5pt}$$
\begin{equation}
                                         \label{02,4,5}
+\sup_{\substack{a\in\bS_{\bar\delta},|b|
\leq K_{0}\\0\leq c\leq K_{0}}}\big[
a_{ij}D_{ij}w^{(\beta)}_{\varepsilon,K}
+b_{i}D_{i}w^{(\beta)}_{\varepsilon,K}
-cw^{(\beta)}_{\varepsilon,K}\big]
\leq\xi_{\varepsilon,K}+\zeta_{\varepsilon,K}\leq0 
\smallskip\end{equation}

\noindent
in $\Pi_{n_{0}-1}$.
 Hence, in $\Pi_{n_{0}-1}$  \smallskip   
$$
\partial_{t}\big(u_{K}^{(\varepsilon)}+w^{\beta}_{\varepsilon,K}\big)+
 H\big[u_{K}^{(\varepsilon)}+w^{\beta}_{\varepsilon,K} \big]
\leq0.
\smallskip\smallskip$$
Since $H[v]-H[v+\gamma \psi]=a^{ij}D_{ij}\psi+b^{i}D_{i}\psi-c\psi$
for some $ \bS_{\delta}$-valued  $a $,  $ \bR^{d}$-valued $b$,
such that $|b|\leq K_{0}$, and $c\geq1$, we have 
in $\Pi_{n_{0}-1}$ that \smallskip
$$
\partial_{t}\big(u_{K}^{(\varepsilon)}+w^{\beta}_{\varepsilon,K}\big)+
 H\big[u_{K}^{(\varepsilon)}+w^{\beta}_{\varepsilon,K}+\beta\psi\big]
\leq-\beta<0.
\smallskip\smallskip$$

This and the definition of $C$-viscosity solutions
imply that,
if $v$
is a continuous in $\bar \Pi$,
$C$-viscosity solution of $\partial_{t}v+H[v]=0$ with boundary data
$g$, then
 the minimum of
$u_{K}^{(\varepsilon)} +w^{ \beta}_{\varepsilon,K}+\beta\psi
-v$ in $\bar\Pi_{n_{0}-1 }$ is either 
nonnegative or is attained
on the parabolic boundary 
of $\Pi_{n_{0}-1}$. 
The same conclusion holds after letting first $\beta\downarrow0$ 
and then
$\varepsilon \downarrow0$. Combining this with the Aleksandrov
estimates showing that $w _{\varepsilon,K}\to0$
as $\varepsilon \downarrow0$ uniformly on $\bar\Pi_{n_{0}}$, we get that
 in $\Pi$  \smallskip
$$
u_{K}-v\geq -\sup_{\Pi\setminus \Pi_{n_{0}-1}}
|u_{K}-v|,
\smallskip$$
which after letting $n_{0}\to\infty$
and then $K\to\infty$
yields $v\leq u$, where $u$ is the common limit
of $u_{K}$, $u_{-K}$, which exists
by Lemma \ref{lemma 2,4,2}. By comparing $v$ with $u_{-K}$, we get
$v\geq u$, and hence uniqueness.
After that estimate \eqref{3,16,6} follows
immediately from \eqref{3,16,7}, and the theorem
is proved.       \qed

\mysection{Proof of Lemma \protect\ref{lemma 2,4,2}}

                                        \label{section 02,1,1}
To prove Lemma  \ref{lemma 2,4,2} we need an auxiliary
result.
In the following theorem   
 $\Omega$ can be just any bounded domain.
Below by $C^{1,2}_{\loc}( \Pi)$ we mean, as usual, the space of
functions $u=u(t,x)$ which are   continuous
in $ \Pi$ along with their derivatives $\partial_{t}u$,
$D_{ij}u$, $D_{i}u$. We recall that $\kappa$ is introduced
before Assumption \ref{assumption 3,13,1} and fix  a constant   
$$
 \tau\in(0,1).
$$

\begin{theorem}
                                               \label{theorem 3.27.1}
Let  
  $u,v\in C^{1,2}_{\loc}( \Pi)\cap C(\bar \Pi)$
be such that for a constant $K\geq1$  \vspace{0pt}
\begin{equation}
                                                 \label{3.30.5}
\partial_{t}u+\max\big(H[u],P[u]-K\big)\geq0\geq
\partial_{t}v+\min\big(H[v],-P[-v]+K\big)
\vspace{10pt}\end{equation}

\noindent
in $\Pi$ and $v\geq u$
on $\partial' \Pi$. Also assume that, for a
constant  $M \in[1,\infty)$,  \vspace{5pt}
\begin{equation}
                                                 \label{3.28.6}
\|u,v\|_{C^{\kappa}(\Pi)}\leq M .
\vspace{5pt}\end{equation}
Then there exist  a  constant    $N \in(0,\infty)$,
depending only on $\tau$,   the diameter of $\Omega$,
    $d$, $K_{0}$, $ \bar H $, and $\delta$,
 and a  constant  $\eta>0$,
depending only on $\tau$,  
    $d$,   and $\delta$,
such that, if $K\geq NM^{ \eta}$
 and
\begin{equation}
                                     \label{5.27.1}
K\geq T^{-1},\quad \Pi^{2/\sqrt{K}}\ne\emptyset,
\end{equation}   
then  in $\Pi$
\begin{equation}
                                                 \label{3.30.1}
u(t,x)-v(t,x)\leq  NK^{-(\kappa-1)/4}
+ NM\omega (M^{-1/\tau}K^{-1} ).
\vspace{5pt}\end{equation}

\end{theorem}

\begin{remark}
                                            \label{remark 4.5.2}
The purpose of introducing $\tau$ is that for $\omega=t^{\tau}$ estimate
\eqref{3.30.1} becomes $u-v\leq  NK^{-(\kappa-1)/4}+NK^{-\tau}$,
which was used  in \cite{Kr15.1} to estimate the rate of convergence    
of finite-difference approximations for \eqref{7.29.10}.
\end{remark}

The statement of this theorem is almost identical
to that of Theorem 3.1 of \cite{Kr15.1} although
that theorem is about the Isaacs equations and our
equations are more general. However, the most part
of the proof follows that of Theorem 3.1 of \cite{Kr15.1}
and, as there, 
we are going to adapt to our situation an argument
from Section 5.A of \cite{CIL}. For that
we need a construction and two lemmas.
 From the start 
throughout the section we will only concentrate
on $K$ satisfying \eqref{5.27.1}.

We take and fix a   function $\zeta=\zeta(t,x)$ as before \eqref{5.3.2}
  and use the
notation $u^{(\varepsilon)}$ introduced in 
\eqref{5.3.2}.
 Recall some standard properties of parabolic mollifiers
in which no regularity properties of $\Omega$ are required:
 If $u\in C^{ \kappa }(\Pi)$, then in $\Pi^{\varepsilon}$
$$
\varepsilon^{ - \kappa } |u-u^{(\varepsilon)} |
+\varepsilon^{ -(\kappa-1)} |Du-Du^{(\varepsilon)} |
\leq N\|u
\|_{C^{  \kappa }(\Pi)},
\vspace{5pt}$$
$$
  |u^{(\varepsilon)}
  |+ |Du^{(\varepsilon)} | +\varepsilon^{2- \kappa }
 |D^{2}u^{(\varepsilon)} |
 +\varepsilon^{2- \kappa }\big|\partial_{t}u^{(\varepsilon)}\big|
\vspace{5pt}$$
\begin{equation}
                                                    \label{3.31.1}
 +\varepsilon^{3- \kappa } |D^{3}u^{(\varepsilon)} |
+\varepsilon^{3- \kappa } \big|D\partial_{t}u^{(\varepsilon)} \big|
+\varepsilon^{4- \kappa } \big|\partial_{t}D^{2}u^{(\varepsilon)} \big|
\leq N\|u
\|_{C^{  \kappa }(\Pi)},
\vspace{5pt}\end{equation}
  where the constants $N$ depend only on $d$ and $\kappa$. 
  
Next, take the constants $\nu,\varepsilon_{0}\in(0,1)$,
specified below in Lemma \ref{lemma 3.30.1} and \eqref{3,17,8},  
respectively, depending only on
    $d$, $K_{0}$, $ \bar H $, $\delta$, and the diameter of $\Omega$,
introduce 
$$
\varepsilon=\varepsilon_{0}M^{-1/(\kappa-1)} K^{-(1-\gamma)/(2\gamma)},
$$
and
consider the function    \vspace{5pt}
$$
W(t,x, y)= u(t,x)-  u^{(\varepsilon)}(t,x)-
\big[  v(t,y)- u^{(\varepsilon)}(t,y)\big]-\nu K
 |x-y|^{2} 
\vspace{5pt}$$
for $ (t,x),(t,y)\in \bar \Pi^{\varepsilon}$.
  Note that $\Pi^{\varepsilon}\ne
\emptyset$ and even $\Pi^{2\varepsilon}\ne\emptyset$
owing to \eqref{5.27.1} and the fact that
$1-\gamma>\gamma$ and $K,M\geq1$.

We will need the following simple observation.
\begin{lemma}
                                                \label{lemma 9,5,1}
For any $\chi>0$ there exists $N=N(\chi,d,\delta)$ such that,
if $K\geq N$, then
\begin{equation}
                                                     \label{9,5,1}
\varepsilon M\leq\chi K^{-(\kappa-1)/4}.
\end{equation}

\end{lemma}

Proof. Since $\varepsilon_{0}<1$, $\kappa\leq2$, and $M\geq1$, the left-hand side
of \eqref{9,5,1} is less than  
$$
K^{-(1-\gamma)/(2\gamma)}=K^{-(5-\kappa)/(14-6\kappa)}.
$$
One easily checks that $(5-\kappa)/(14-6\kappa)>1/2>(\kappa-1)/4$ for  
$\kappa\in(1,2)$ and this proves the lemma.    \qed

Denote by  $(\bar t, \bar x,  \bar y)$ a maximum point
of $W$ in $ [0,T-
\varepsilon^{2}] \times (\bar \Omega^{\varepsilon})^{2}$.
Observe that, obviously,  \vspace{5pt}
$$
  u(\bar t,\bar x)-  u^{(\varepsilon)}(\bar t,\bar x)-
\big[  v(\bar t,\bar y)-  u^{(\varepsilon)}(\bar t,\bar y)\big]-\nu K
 |\bar x-\bar y|^{2} 
\vspace{5pt}$$
$$
\geq
  u(\bar t,\bar x)-  u^{(\varepsilon)}(\bar t,\bar x)-
\big[ v(\bar t,\bar x)-  u^{(\varepsilon)}(\bar t,\bar x)\big],
\vspace{5pt}$$
which implies that
\begin{equation}
                                                    \label{4.21.2}
 |\bar x-\bar y|\leq N M/(\nu K) .
\vspace{5pt}\end{equation}
where and below by $N$  with indices or without them we denote 
various constants
depending only on 
    $d$, $K_{0}$, $ \bar H $, $\delta$, $\tau$, and  the diameter of $\Omega$,
unless specifically stated otherwise.
By the way, recall that $ \kappa $ and, hence, $\gamma$
 depend  only on $d$ and $\delta$.

\begin{lemma}
                                              \label{lemma 3.30.1}
There exist  a constant   $\nu\in(0,1)$,
depending only on
    $d$, $K_{0}$, $ \bar H $, $\delta$, and the diameter of $\Omega$, and 
a constant $N\in[1,\infty)$
  such that, if  \vspace{5pt}
\begin{equation}
                                            \label{3.31.09}
K\geq N\nu^{-\eta_{1}}\varepsilon_{0}^{( \kappa-2 ) \eta_{1}}
M^{ \eta_{1}(2-\kappa)/(\kappa-1)},
\vspace{5pt}\end{equation}
where $  \eta_{1}^{-1}=1-(2- \kappa )(1-\gamma) (2\gamma)^{-1}$ \($>0$\),
and $ \bar x ,\bar y \in \Omega^{\varepsilon}$ 
and $ \bar t<T-\varepsilon^{2}$, then

\(i\,\) we have \vspace{5pt}
\begin{equation}
                                            \label{3.27.6}
2\nu  |\bar x-\bar y|\leq NMK^{-1} \varepsilon^{ \kappa-1 }
=N\varepsilon^{ \kappa-1 }_{0}K^{-(3+\kappa)/(8\gamma)},
\quad|\bar x-\bar y|\leq\varepsilon/2;
\vspace{5pt}\end{equation}

\(ii\,\) for any $\xi,\eta\in\bR^{d}$ \vspace{5pt}
\begin{equation}
                                            \label{3.27.7}
D_{ij} [  u -  u^{(\varepsilon)}  ](\bar t,\bar x)\xi^{i}\xi^{j}
-D_{ij} [  v - u^{(\varepsilon)}  ](\bar t,\bar y)\eta^{i}\eta^{j}
\leq 2\nu K|\xi-\eta|^{2},
\vspace{5pt}\end{equation}
\begin{equation}
                                            \label{3.31.2}
\partial_{t} [  u - u^{(\varepsilon)}  ]
(\bar t,\bar x) \leq 
 \partial_{t} [ v -  u^{(\varepsilon)}  ]
(\bar t,\bar y) ;
\vspace{5pt}\end{equation}

\(iii\,\) we have
\begin{equation}
                                               \label{3.27.8}
\partial_{t}  u(\bar t,\bar x)+
H[u](\bar t,\bar x)\geq0,
\vspace{10pt}\end{equation} 
\begin{equation}
                                               \label{3.28.1}
\partial_{t}  v(\bar t,\bar y)+H[v]
(\bar t,\bar y)\leq0.
\vspace{5pt}\end{equation}

\end{lemma}

Proof. The first inequality in
\eqref{3.27.6}  follows from \eqref{3.31.1} and
the fact that the first derivatives
of $W$ with respect to $x$ vanish  at $\bar x$, that is,   \vspace{5pt}
$$
D [\bar u-\bar u^{(\varepsilon)}
](\bar t,\bar x)=2\nu K(\bar x-\bar y).
\vspace{5pt}$$ 
Also the matrix of second-order derivatives of $W$
with respect to $(x,y)$
is nonpositive at $( \bar t,\bar x,\bar y)$ as well as its 
(at least one sided if $\bar t=0$)
derivative
with respect to   $t$, which yields (ii).

By taking $\eta=0$ in \eqref{3.27.7} and using the fact that
$ |D^{2} u^{(\varepsilon)} |\leq N M
\varepsilon^{ \kappa-2}$ we see that  \vspace{5pt}
$$
D^{2} u(\bar t,\bar x)\leq 2\nu K+NM
\varepsilon^{ \kappa-2}.
$$
 
Similarly,  \vspace{5pt}
$$
D^{2}  v(\bar t,\bar x)\geq -N (\nu K+ M
\varepsilon^{ \kappa-2} ) ,
\vspace{5pt}$$
which yields  \vspace{5pt}
$$
P[u](\bar t,\bar x)
\leq N_{1} (\nu K+ M
\varepsilon^{\kappa-2}  ),\quad -P[-v](\bar t,\bar y)
\geq -
 N_{1} (\nu K+ M
\varepsilon^{\kappa-2} ).
\vspace{5pt}\vspace{5pt}$$

  Since
 $
H(\sfu,t,x)=a^{ij}\sfu_{ij}''+H(\sfu',0,t,x)$
where $(a^{ij})\in\bS_{\delta}$ ($(a^{ij})$ depends on $\sfu,t,x$)
and $|H(\sfu',0,t,x)|\leq K_{0}|\sfu'|
+\bar H$ and $M\geq1$ and $\varepsilon<1$, we also have (by increasing
 the above $N_{1}$ if necessary)    \smallskip
\begin{equation}
                                                 \label{5.19.1}
H[u](\bar t,\bar x)
\leq N_{1} (\nu K+ M
\varepsilon^{\kappa-2}  ),\quad H[v](\bar t,\bar y)
\geq -
 N_{1} (\nu K+ M
\varepsilon^{\kappa-2} ).
\vspace{5pt}\end{equation}
Also it follows from \eqref{3.31.2} and \eqref{3.31.1} that \vspace{5pt}
\begin{equation}
                                                 \label{5.19.2}
\partial_{t}u(\bar t,\bar x)\leq \partial_{t}v(\bar t,\bar y)
+N_{2} M
\varepsilon^{\kappa-2} .
\vspace{5pt}\end{equation}

Now, if $H[u](\bar t,\bar x)\leq P[u](\bar t,\bar x)-K$,
then at $(\bar t,\bar x)$    \vspace{5pt}
$$
0\leq
\partial_{t}u 
+\max\big(H[u] , P[u] -K\big)\leq\partial_{t}u+
N_{1} (\nu K+ M
\varepsilon^{\kappa-2}  )-K,
\vspace{5pt}$$
$$
\partial_{t}u \geq K-N_{1} (\nu K+ M
\varepsilon^{\kappa-2} )
\vspace{5pt}$$
and at $(\bar t,\bar y)$   \vspace{5pt}
$$
0\geq\partial_{t}v 
+\min\big(H[v] , -P[-v]+K\big)
\geq K-2N_{1} (\nu K+ M
\varepsilon^{\kappa-2}  )-N_{2} M
\varepsilon^{\kappa-2}.
\vspace{5pt}$$
Here $M\varepsilon^{\kappa-2}\leq \nu K$, which is equivalent to 
\eqref{3.31.09} with $N=1$. Hence \smallskip
$$
K\leq 4N_{1}\nu K+N_{2}\nu K,
\medskip$$
 which is impossible if we
  choose and fix $\nu$  such that  \vspace{5pt}
$$
(4N_{1} +N_{2})\nu\leq1/2.
\vspace{5pt}$$

It follows that $H[u](\bar t,\bar x)\geq P[u](\bar t,\bar x)-K$, \vspace{5pt}
$$
\partial_{t}u(\bar t,\bar x)+H[u](\bar t,\bar x)\geq0,
\vspace{5pt}$$
and this proves \eqref{3.27.8}.

Similarly, if $-P[-v](\bar t,\bar y)+K\leq H[v](\bar t,\bar y)$,
then at $(\bar t,\bar y)$   \vspace{5pt}
$$
0\geq\partial_{t} v -
 N_{1} (\nu K+ M
\varepsilon^{\kappa-2}  )+K,
\vspace{5pt}$$
and at $(\bar t,\bar x)$ we have $\partial_{t} u\leq
N_{1} (\nu K+ M
\varepsilon^{\kappa-2}  ) -K+N_{2} M
\varepsilon^{\kappa-2}$ and  \vspace{5pt}
$$
0\leq 
\partial_{t} u+\max\big(H[u],P[u]-K\big)\leq
-K+ 2N_{1} (\nu K+ M
\varepsilon^{\kappa-2}  )+N_{2} M
\varepsilon^{\kappa-2},
\vspace{5pt}$$
which again is impossible with the above choice of $\nu$
for $K$ satisfying \eqref{3.31.09}. This yields
\eqref{3.28.1}.

Moreover, not only $M\varepsilon^{\kappa-2}\leq \nu K$ for
$K$
satisfying \eqref{3.31.09}, but we also have
$NM\varepsilon^{\kappa-2}\leq \nu K$, where $N$ is taken from
\eqref{3.27.6}, if we increase $N$ in \eqref{3.31.09}.
This yields the second inequality in \eqref{3.27.6}.
 
 The lemma is proved.   \qed

Everywhere below in this section $\nu$
is the constant from Lemma \ref{lemma 3.30.1}.

\begin{lemma}
                                            \label{lemma 4.23.1}
There exist   a   constant $N$ 
such that,
for any $\mu\geq0$,
  if $K\geq NM^{ \eta_{2}}$, where $\eta_{2}=8/(5-\kappa)$, and     \vspace{5pt}
\begin{equation}
                                                 \label{3.28.5}
W(\bar t,\bar x, \bar y)\geq 2K^{-(\kappa-1)/4}
+\mu M\omega (M^{-1/\tau}K^{-1} ),
\end{equation}
then 
\begin{equation}
                                            \label{3.27.4}
 u(\bar t,\bar x)- 
   v(\bar t,\bar y) -
\nu K|\bar x-\bar y|^{2} 
\geq K^{-(\kappa-1)/4}+\mu M\omega (M^{-1/\tau}K^{-1} ).
\vspace{5pt}\end{equation}
Furthermore,
$ \bar x,\bar y \in \Omega^{2\varepsilon}$ and
$\bar t <T-\varepsilon^{2}$.
\end{lemma}

Proof.  
It follows from \eqref{4.21.2} 
that  (recall that $\nu$ is already fixed)  \vspace{5pt}
 \begin{equation}
                                                 \label{4.21.5}
\big|  u^{(\varepsilon)}(\bar t,\bar x)-
 u^{(\varepsilon)}(\bar t,\bar y)\big|\leq M  |\bar x-\bar y|
\leq N_{0} M^{2}/  K.
\vspace{5pt}\end{equation}
Hence we have from \eqref{3.28.5} that    \vspace{5pt}
$$
 u(\bar t,\bar x)- 
   v(\bar t,\bar y) -
\nu K|\bar x-\bar y|^{2}\geq 2 K^{-(\kappa-1)/4}-N_{0} M^{2}/  K
+\mu M\omega (M^{-1/\tau}K^{-1} ),
\vspace{5pt}$$
and \eqref{3.27.4} follows provided that    \vspace{5pt}
\begin{equation}
                                                   \label{4.22.2}
N_{0}M^{2}/K\leq (1/4) K^{-(\kappa-1)/4},
\vspace{5pt}\end{equation}
which indeed holds if
\begin{equation}
                                                   \label{3.31.6}
K\geq N 
M^{ \eta_{2}}.
\vspace{5pt}\end{equation}

Now, if $\bar x$ or $\bar y$ is 
 in $\bar \Omega^{\varepsilon}
\setminus \Omega^{2\varepsilon}$, then
for appropriate $\hat x\in\partial \Omega$ and $\hat y\in\partial \Omega $
either 
$$
  u(\bar t,\bar x)- 
   v(\bar t,\bar y)\leq u(\bar t,\bar x)- 
   u(\bar t,\hat x)+
  v(\bar t,\hat x)-  v(\bar t,\bar y)
\leq M \big(4\varepsilon+|\bar x-\bar y|\big)
$$
or
$$
  u(\bar t,\bar x)- 
   v(\bar t,\bar y)\leq \bar u(\bar t,\bar x)
-  u(\bar t,\hat y)
+   v(\bar t,\hat y)-  v(\bar t,\bar y)
\leq M \big(4\varepsilon+|\bar x-\bar y| \big).
\vspace{5pt}$$

In any case in light of \eqref{3.27.4}, \eqref{4.21.2},
and \eqref{4.22.2} \vspace{5pt}
$$
4\varepsilon M +N_{0} M^{2}/K-
\nu K|\bar x-\bar y|^{2} \geq K^{-(\kappa-1)/4},
$$
\begin{equation}
                                                   \label{4.23.5}
4\varepsilon M   \geq (3/4)K^{-(\kappa-1)/4}.
\vspace{5pt}\end{equation}

  Owing to Lemma \ref{lemma 9,5,1},
  inequality \eqref{4.23.5} is 
impossible if $K\geq N$ or, upon adjusting $N$ in \eqref{3.31.6}
appropriately, if \eqref{3.31.6} holds (recall that $M\ge1$).
Below we assume that \eqref{3.31.6} holds after the modification,
so that 
we have  $ \bar x,\bar y \in \Omega^{\varepsilon}$. 

Furthermore, if $\bar t= T-\varepsilon^{2}$, then
(recall \eqref{4.21.5} and that $  u\leq v$
on $\partial' \Pi$)
$$
W(\bar t,\bar x, \bar y)\leq 
NM\varepsilon^{ \kappa }+  v(\bar t,\bar x)-
 v(\bar t,\bar y)+NM^{2}/K
\vspace{5pt} $$
$$
\leq NM\varepsilon+M |\bar x-\bar y| 
 +NM^{2}/K\leq NM^{2}/K+NM\varepsilon ,
\vspace{5pt}\vspace{5pt}$$
which is less than $(1/2)K^{-(\kappa-1)/4}$ in light of \eqref{4.22.2}
and    Lemma \ref{lemma 9,5,1}, 
  again after 
perhaps further adjusting the constant in \eqref{3.31.6}.
This is impossible due to \eqref{3.28.5}.
Hence, $\bar t <T-\varepsilon^{2}$ and this
finishes the proof of the lemma.      \qed

We also need a simple result based on solving quadratic inequalities.

\begin{lemma}
                                              \label{lemma 3,17,2}
 Let $\beta\in[0,\delta/2]$. Then for any $\alpha\geq-4\nu$ we have  \vspace{5pt}
$$
\beta\alpha-\delta\frac{\alpha^{2}}{\alpha+4\nu}\leq
2(\nu/\delta)\beta^{2}.
\vspace{5pt}$$

\end{lemma}

{\bf Proof of Theorem \ref{theorem 3.27.1}}.
Fix a (large) constant $\mu>0$ to be specified later
as a constant, depending only on  $d$, $K_{0}$, 
 $\bar H$, $\delta$, $\tau$,
and the diameter of $\Omega$,
recall that $\nu$ is found in Lemma \ref{lemma 3.30.1} and
first assume that  \vspace{5pt}
$$
W(t,x ,y)\leq 2K^{-(\kappa-1)/4}+\mu M\omega (M^{-1/\tau}K^{-1} )
\vspace{5pt}$$
for  $ (t,x),(t,y)\in \bar \Pi^{\varepsilon}$. Observe that for 
any point $(t,x)\in \bar \Pi $
one can find a point $(s,y)\in \bar \Pi^{\varepsilon}$
with $|x-y|\leq\varepsilon$ and $|t-s|\leq\varepsilon^{2}$ and
then   \vspace{5pt}
$$
  u(t,x)- v(t,x)\leq   u(s,y)-  v(s,y)
+NM \varepsilon
\vspace{5pt} $$
$$
\leq W(s,y, y)+NM \varepsilon
\leq 2K^{-(\kappa-1)/4}+NM \varepsilon
+\mu M\omega (M^{-1/\tau}K^{-1} ).
\vspace{5pt}$$
In that case,  as follows from  Lemma \ref{lemma 9,5,1}, \eqref{3.30.1} holds for 
$K$ satisfying \eqref{3.31.6} with any $N\geq1$ in \eqref{3.31.6}.

It is clear now that, to prove the theorem, it suffices to find  
  $\mu$ such that the inequality  \eqref{3.28.5}
is impossible  
if $K\geq NM^{ \eta}$ with $N$
and $\eta$ as in the statement of the theorem
and at least   not smaller than those in \eqref{3.31.6}.
 Of course, we will argue by 
contradiction and suppose that \eqref{3.28.5} holds,
despite \eqref{3.31.6} is valid and
\eqref{3.31.09} is satisfied with $\nu$ fixed in Lemma \ref{lemma 3.30.1}
and $\varepsilon_{0}\in(0,1)$, which is yet to be specified.

Then 
\eqref{3.27.4} holds,
and, in particular,
\begin{equation}
                                                     \label{3,16,4}
u(\bar t,\bar x)\geq v(\bar t,\bar y).
\end{equation}
Also by Lemma \ref{lemma 4.23.1} the points
   $ \bar x ,\bar y$ are in $  \Omega^{\varepsilon}$ 
(even in $\Omega^{2\varepsilon}$)
and $\bar t<T-\varepsilon^{2}$, so that
we can use   the conclusions of Lemma \ref{lemma 3.30.1}.

Denote 
$$
A=D^{2} [u-u^{(\varepsilon)} ](\bar t,\bar x),\quad B=
D^{2} [v-u^{(\varepsilon)} ](\bar t,\bar y)
\vspace{5pt}$$
and
interpret matrices as linear operators in a usual way and constants
as operators of multiplications by these constants.
Observe that \eqref{3.27.7} implies that
the operator $B+2\nu K$ is nonnegative and, hence,
$B+4\nu K$ is strictly positive. Then for
$$
\eta=4\nu K(B+4\nu K)^{-1}\xi,\quad \eta-\xi=-B(B+4\nu K)^{-1}\xi
\vspace{5pt}$$
inequality
\eqref{3.27.7} yields  \vspace{5pt}
\begin{align*}
\langle A\xi,\xi\rangle\leq &\,
 \big\langle B4\nu K(B+4\nu K)^{-1}\xi,
4\nu K(B+4\nu K)^{-1}\xi\big\rangle
\\[10pt]
&\,+4\nu K\big|B(B+4\nu K)^{-1}\xi\big|^{2}
\\[10pt]
=&\,\big\langle 4\nu KB(B+4\nu K)^{-1}\xi,\xi\big\rangle.
 \end{align*}
Hence $A\leq 4\nu KB(B+4\nu K)^{-1}$ and  \vspace{5pt}
\begin{align*}
D^{2} u(\bar t,\bar x)\leq &\, D^{2}u^{(\varepsilon)} (\bar t,\bar x)
+4\nu KB(B+4\nu K)^{-1}
\\[10pt]
=  &\, D^{2}v  (\bar t,\bar y)-B+
D^{2}u^{(\varepsilon)} (\bar t,\bar x)
-D^{2}u^{(\varepsilon)} (\bar t,\bar y)+4\nu KB(B+4\nu K)^{-1} 
\\[10pt]
= &\, D^{2}v  (\bar t,\bar y)-B^{2}(B+4\nu K)^{-1}+D^{2}u^{(\varepsilon)} (\bar t,\bar x)
-D^{2}u^{(\varepsilon)} (\bar t,\bar y).
\vspace{5pt}\end{align*}
 
We now use  \eqref{3.27.6} to get that
$
 |\bar x-\bar y| 
\leq  NMK^{-1}
 \varepsilon^{ \kappa-1 } 
$,
and in light of \eqref{3.31.1} that
$$
\big|D^{2} u^{(\varepsilon)}(\bar t,\bar x)-
D^{2} u^{(\varepsilon)}(\bar t,\bar y)\big|\leq 
NM\varepsilon^{ \kappa-3}|\bar x-\bar y|\leq NC,
$$
where
$$
C=M^{2}K^{-1}\varepsilon^{2\kappa-4}.
$$

Also \eqref{3.31.2} reads   \vspace{5pt}
$$
\partial_{t}  u (\bar t,\bar x)
\leq\partial_{t}v (\bar t,\bar y)
+\partial_{t}  u^{(\varepsilon)}(\bar t,\bar x)
-\partial_{t}
 u^{(\varepsilon)}(\bar t,\bar y)
\vspace{5pt}$$
and
as is easy to see
$$
\big|\partial_{t} u^{(\varepsilon)}(\bar t,\bar x)
-\partial_{t}
 u^{(\varepsilon)}(\bar t,\bar y)\big|\leq  
NC.
$$

Thus far, we have
$$
D^{2} u(\bar t,\bar x)
\leq D^{2}v  (\bar t,\bar y)-B^{2}(B+4\nu K)^{-1}+NC,
\vspace{5pt}\vspace{5pt}$$
\begin{equation}
                                              \label{3,16,2}
 \partial_{t}  u(\bar t,\bar x)
\leq \partial_{t} v  (\bar t,\bar y) +NC.
\end{equation}

Next,
$$
D_{i}  u(\bar t,\bar x)=2\nu K(\bar x^{i}-\bar y^{i})+
D_{i} u^{(\varepsilon)}(\bar t,\bar x),
\vspace{5pt}$$
$$
 D_{i}  v(\bar t,\bar y)=
2\nu K(\bar x^{i}-\bar y^{i})
+D_{i}  u^{(\varepsilon)}(\bar t,\bar y),
\vspace{5pt}$$
$$
D_{i}  u(\bar t,\bar x)-D_{i}\bar v(\bar t,\bar y)=
D_{i}  u^{(\varepsilon)}(\bar t,\bar x)-
D_{i}  u^{(\varepsilon)}(\bar t,\bar y),
\vspace{5pt}$$
where
$$
\big|D   u^{(\varepsilon)}(\bar t,\bar x)
-D\ u^{(\varepsilon)}(\bar t,\bar y)\big|
\leq NM \varepsilon^{ \kappa-2}|\bar x-\bar y| \leq
NC
\vspace{5pt}$$
and, therefore,
\begin{equation}
                                              \label{3,16,3}
\big|D  u(\bar t,\bar x)-D  v(\bar t,\bar y)\big|\leq NC.
\vspace{5pt}\end{equation}

It follows from \eqref{3,16,2}, \eqref{3,16,3},  
\eqref{3.27.8}, and \eqref{3,16,4} by Assumption \ref{assumption 3,13,1}
and the fact that $D_{\sfu'_{0}}H\leq-1$ (see  
\eqref{3,15,5}) that
\begin{align*}
0\leq &\,\partial_{t}u(\bar t,\bar x)
+H\big(u(\bar t,\bar x),Du(\bar t,\bar x),
D^{2}v  (\bar t,\bar y)-B^{2}(B+4\nu K)^{-1}
 +NC,\bar t,\bar x\big)
\\[10pt]
\leq &\,\partial_{t}v (\bar t,\bar y)
+NC+H
\big(v(\bar t,\bar y),Dv(\bar t,\bar y),
D^{2}v  (\bar t,\bar y)
 +NC,\bar t,\bar x\big)
\\[10pt]
&\,-\big[u(\bar t,\bar x)-v(\bar t,\bar y)\big]-\delta\tr  B^{2}(B+4\nu K)^{-1} 
\\[10pt]
\leq &\,\partial_{t}v (\bar t,\bar y)
+NC+H
\big(v(\bar t,\bar y),Dv(\bar t,\bar y),
D^{2}v  (\bar t,\bar y),\bar t,\bar x\big)
\\[10pt]
&\,-\big[u(\bar t,\bar x)-v(\bar t,\bar y)\big]-\delta\tr  B^{2}(B+4\nu K)^{-1}
+NC.
\vspace{5pt}\end{align*}
This along with \eqref{3.28.1} and Assumption 
\ref{assumption 3,13,1} (iii) yields (recall that $M\geq1$)   \vspace{5pt}
\begin{align*}
u(\bar t,\bar x)-v(\bar t,\bar y)
\leq &\, K_{0}|\bar x-\bar y|^{\gamma}\big |D^{2}v  (\bar t,\bar y)\big |
+NM\omega\big(|\bar x-\bar y|\big)
\\[10pt]
 &\,-\delta\tr  B^{2}(B+4\nu K)^{-1}
+NC.
\vspace{5pt}\end{align*}

Upon combining this with \eqref{3.27.4} we arrive at
\begin{align}
 2K^{-(\kappa-1)/4}+\mu M\omega\big(M^{-1/\tau}K^{-1}\big)
\leq K_{0}|\bar x-\bar y|^{\gamma}\big |D^{2}v  (\bar t,\bar y)\big |+NC
\nonumber\\[10pt]
                                              \label{3,16,5}
+N_{1}M\omega\big(|\bar x-\bar y|\big)-\delta\tr  B^{2}(B+4\nu K)^{-1}.
\vspace{5pt}\end{align}
 
Now we choose $\mu=N_{1}$ and observe that (see \eqref{3.27.6}) \vspace{5pt}
$$
 |\bar x-\bar y|\leq N_{2}
K^{-(3+\kappa)/(8\gamma)}\leq M^{-1/\tau}K^{-1}
 $$
for
\begin{equation}
                                              \label{3,17,1}
K^{\theta_{1}}\geq N_{2}M^{1/\tau},
\vspace{5pt}\end{equation}
where $\theta_{1}=(\kappa-1)(5-\kappa)(14-6\kappa)^{-1}$.
Then for such $K$ we infer from \eqref{3,16,5} that  \vspace{5pt}
\begin{equation}
                                              \label{3,17,2}
2K^{-(\kappa-1)/4} 
\leq K_{0}|\bar x-\bar y|^{\gamma}  \big|D^{2}v  (\bar t,\bar y) \big |
 -\delta\tr  B^{2}(B+4\nu K)^{-1}
+N_{3}C.
\vspace{5pt}\end{equation}

Here, owing to \eqref{3.31.1} and \eqref{3.27.6}, \smallskip
$$
|\bar x-\bar y|^{\gamma}  \big|D^{2}v  (\bar t,\bar y)  \big|   
\leq N\varepsilon^{ (\kappa-1)\gamma }_{0}K^{-(3+\kappa)/8}\big(|B|
+\big|D^{2}u^{(\varepsilon)}(\bar t,\bar y)  \big|\big)
$$
\medskip$$
\leq N\varepsilon^{ (\kappa-1)\gamma }_{0}K^{-(3+\kappa)/8} |B|
+N_{4}M^{1/(\kappa-1)}\varepsilon_{0}^{ (\kappa-1)\gamma +\kappa-2}
K^{\theta_{2} (\kappa-1)/4 },
\medskip$$
where $\theta_{2} =(7\kappa-19)/(14-6\kappa)$, and as is easy to see
$\theta_{2} \leq-3/2$ for 
 $\kappa\in[1,2]$, so that
$$
N_{4}M^{1/(\kappa-1)}\varepsilon_{0}^{ (\kappa-1)\gamma +\kappa-2}
K^{\theta_{2}(\kappa-1)/4 }\leq  (1/4) K^{-(\kappa-1)/4}
$$
if, for instance,
\begin{equation}
                                              \label{3,17,4}
K \geq \big(4N_{4}M^{1/(\kappa-1)}
\varepsilon_{0}^{ (\kappa-1)\gamma +\kappa-2}\big)^{8/(\kappa-1)}.
\end{equation}

In what concerns the last term in \eqref{3,17,2}, note that \smallskip
$$
N_{3}C =N_{3}M^{2/(\kappa-1)}\varepsilon^{2\kappa-4}_{0}K^{\theta_{3} (\kappa-1)/4},
\medskip$$
where $\theta_{3} =(4 \kappa-12)/(7-3\kappa)\leq-2$  for 
 $\kappa\in[1,2]$, so that   \smallskip
$$
N_{3}C \leq (1/4) K^{-(\kappa-1)/4}
$$
if
\begin{equation}
                                              \label{3,17,5}
K\geq\big(4N_{3}M^{2/(\kappa-1)}\varepsilon^{2\kappa-4}_{0}\big)^{4/(\kappa-1)}.
\vspace{5pt}\vspace{5pt}\end{equation}

We conclude that, for $K$ satisfying
   \eqref{3,17,4} 
and \eqref{3,17,5}, relation \eqref{3,17,2} yields  \vspace{5pt}
\begin{equation}
                                              \label{3,17,7}
K^{-(\kappa-1)/4} \leq
N  
\varepsilon^{ (\kappa-1)\gamma}_{0}K^{-(3+\kappa)/8} |B |
 -\delta\tr  B^{2}(B+4\nu K)^{-1}.
\vspace{5pt}\end{equation}

Next, observe that, by Lemma \ref{lemma 3,17,2} applied
after we diagonalize $B$ and set $B=K\alpha$ implies that
the right-hand side of \eqref{3,17,7} is   \vspace{5pt}
$$
K\Big[N  
\varepsilon^{ (\kappa-1)\gamma}_{0}K^{-(3+\kappa)/8} |\alpha |
 -\delta\tr  \alpha^{2}(\alpha+4\nu  )^{-1}\Big]
\vspace{5pt}$$
$$
\leq N_{5} \varepsilon_{0}^{2(\kappa-1)\gamma}K^{-(\kappa-1)/4},
\medskip$$
where the inequality holds owing to Cauchy's inequality.
We can certainly assume that $N_{5}\geq1$
and then 
  we can choose $\varepsilon_{0}\in(0,1)$ so that \vspace{5pt}
\begin{equation}
                                                  \label{3,17,8}
N_{5} \varepsilon_{0}^{2(\kappa-1)\gamma}=1/2,
\vspace{5pt}\end{equation}
which along with \eqref{3,17,7} leads to the desired 
contradiction:  $K^{-(\kappa-1)/4}\leq(1/2)K^{-(\kappa-1)/4}$.

With so specified $\varepsilon_{0}$  we rewrite   
condition \eqref{3.31.09} (with $\nu$ fixed in Lemma \ref{lemma 3.30.1}),
  conditions \eqref{3.31.6}, \eqref{3,17,1}, 
\eqref{3,17,4}, and \eqref{3,17,5}  as
$$
K\geq N M^{ \eta_{1}(2-\kappa)/(\kappa-1)},\quad K\geq N 
M^{ \eta_{2}},\quad K \geq N  M^{1/(\theta_{1}\tau)},
$$
$$
K\geq N M^{ 8/(\kappa-1)^{2}} .
$$
 
Since $M\geq1$, for $\eta'$ defined as the sum of the 
above powers of $M$ and $N'$ defined as the sum of the 
above $N$'s, the inequality  \eqref{3.28.5}
is impossible for
  $K\geq N'M^{\eta'}$ and this brings the proof of the theorem     
to an end.         \qed

{\bf Proof of Lemma \ref{lemma 2,4,2}}. 
Let $\Omega_{n}$, $n=2,3,...$, be a sequence of 
strictly expanding smooth
domains whose union is $\Omega$ and set $\Pi_{n}=[0, T(1-1/n))
\times\Omega_{n}$. Then for any $n_{0}=3,4,...$ and all sufficiently small
$\varepsilon>0$
$$
\xi_{\varepsilon, K} =\partial_{t}u_{  K}^{(\varepsilon)}
+\max \big(H\big[u_{ K}^{(\varepsilon)}\big],
 P\big[u_{  K}^{(\varepsilon)}\big]-K \big), 
\vspace{5pt}\vspace{5pt}$$
$$
\xi_{\varepsilon,- K} =\partial_{t}u_{- K}^{(\varepsilon)}
+\min\big (H\big[u_{-K}^{(\varepsilon)}\big],
- P\big[-u_{-K}^{(\varepsilon)}\big]+ K \big), 
\vspace{5pt}$$
  are well defined in $\Pi_{n_{0}}$. Since the second-order derivatives  
with respect to $x$ and the first derivative 
with respect to $t$
of $u_{ \pm K}$    are in $L_{p}(\Pi^{\rho})$  
for any $p$ and $\rho$,
we have $ \xi_{\varepsilon, \pm K} 
\to0$ as $\varepsilon\downarrow0$
in any $L_{p}(\Pi_{n_{0}})$ for any $K$
and any $p>1$. Furthermore,
$ \xi_{\varepsilon,\pm K} $ are
 continuous {\em because $H(\sfu,t,x)$
is continuous\/}. Therefore, there exist smooth functions  
$ \zeta_{\varepsilon,K} $ on $\bar \Pi_{n_{0}}$ such that
$$
-\varepsilon\leq
  \zeta_{\varepsilon,K}-
\min\big(\xi_{\varepsilon,K},-\xi_{\varepsilon,-K}\big)\leq 0
$$
in $\Pi_{n_{0}}$ for all small $\varepsilon>0$.

Since $\Omega_{n_{0}}$ is smooth, by Theorem 1.1
of \cite{DKL_12} there exists a unique
 $
w_{\varepsilon,K}\in \bigcap_{ 
p>1}W^{1, 2}_{p }(\Pi_{n_{0}} ) 
 $
 satisfying 
$$
\partial_{t}w_{\varepsilon,K}+
\sup_{\substack{a\in\bS_{\bar\delta},|b|\leq K_{0}\\0\leq
 c \leq K_{0}}}\big[
a_{ij}D_{ij}w_{\varepsilon,K}+b_{i}D_{i}w_{\varepsilon,K}    
-cw_{\varepsilon,K}\big]=
\zeta_{\varepsilon,K}
\smallskip\smallskip$$
 in $\Pi_{n_{0}} $  \(a.e.\)  and such
that $w_{\varepsilon,K}=0$ on $\partial'\Pi_{n_{0}} $. By the maximum principle
such $w_{\varepsilon,K}$ is unique.     Then owing to the continuity
of $\zeta_{\varepsilon,K}$, for any $\varepsilon$ and $K$
for all sufficiently small $\beta>0$, we have   \smallskip
$$
\partial_{t}w^{(\beta)}_{\varepsilon,K}+
\sup_{\substack{a\in\bS_{\bar\delta},|b|\leq K_{0}\\0\leq c
\leq K_{0}}}\big[
a_{ij}D_{ij}w^{(\beta)}_{\varepsilon,K}+
b_{i}D_{i}w^{(\beta)}_{\varepsilon,K}
-cw^{(\beta)}_{\varepsilon,K}\big]
\leq
\zeta^{(\beta)}_{\varepsilon,K}\leq \zeta_{\varepsilon,K} +\varepsilon
$$
in $\Pi_{n_{0}-1}$.

 For $w^{\beta}_{\varepsilon,K}(t,x):=w^{(\beta)}_{\varepsilon,K}
+\varepsilon(T-t)$ now similarly to \eqref{02,4,5}   \smallskip
$$
\partial_{t}\big(u_{K}^{(\varepsilon)}- w^{ \beta }_{\varepsilon,K}\big)+
\max\big(H\big[u_{K}^{(\varepsilon)}-w^{ \beta }_{\varepsilon,K}\big],
P\big[u_{K}^{(\varepsilon)}-w^{ \beta }_{\varepsilon,K}\big]-K\big) 
\vspace{5pt}\vspace{5pt}$$
$$
\geq \partial_{t}u_{K}^{(\varepsilon)}
+\max\big(H\big[u_{K}^{(\varepsilon)}\big],
P\big[u_{K}^{(\varepsilon)}\big]-K\big)+
\varepsilon
\vspace{5pt}\vspace{5pt}$$
$$
-
\sup_{\substack{a\in\bS_{\bar\delta},|b|
\leq K_{0}\\0\leq c\leq K_{0}}}\big[
a_{ij}D_{ij}w^{(\beta)}_{\varepsilon,K}+b_{i}D_{i}
w^{(\beta)}_{\varepsilon,K}-cw^{(\beta)}_{\varepsilon,K}\big]
\geq\xi_{\varepsilon,K}-\zeta_{\varepsilon,K}\geq0,
\vspace{5pt}\vspace{5pt}$$
$$
\partial_{t}\big(u_{-K}^{(\varepsilon)}
+w^{ \beta}_{\varepsilon,K}\big)+
\min\big(H\big[u_{-K}^{(\varepsilon)}+w^{ \beta }_{\varepsilon,K}\big],
-P\big[-u_{-K}^{(\varepsilon)}-w^{ \beta }_{\varepsilon,K}\big]+K\big) 
\vspace{5pt}\vspace{5pt}$$
$$
\leq \partial_{t}u_{-K}^{(\varepsilon)}+
\min\big(H\big[u_{-K}^{(\varepsilon)}\big],
-P\big[-u_{-K}^{(\varepsilon)}\big]+K
\big)+\partial_{t}w^{ \beta }_{\varepsilon,K}
\vspace{5pt}\vspace{5pt}$$
$$
+\sup_{\substack{a\in\bS_{\bar\delta},|b|
\leq K_{0}\\0\leq c\leq K_{0}}}\big[
a_{ij}D_{ij}w^{(\beta)}_{\varepsilon,K}+
b_{i}D_{i}w^{(\beta)}_{\varepsilon,K}
-cw^{(\beta)}_{\varepsilon,K}\big]
\leq\xi_{\varepsilon,-K}+\zeta_{\varepsilon,K}\leq 0
\vspace{5pt}$$
in $\Pi_{n_{0}-1}$.

After setting
$$
\mu^{ \beta}_{\varepsilon,K}=\sup_{\partial'
\Pi_{n_{0}-1}}\big(
u_{K}^{(\varepsilon)}-u_{-K}^{(\varepsilon)}
-2w^{ \beta}_{\varepsilon,K}\big)_{+}
$$
we conclude by Theorem \ref{theorem 3.27.1} applied to
$u_{K}^{(\varepsilon)}-w^{ \beta}_{\varepsilon,K}$
and $u_{-K}^{(\varepsilon)}+w^{ \beta}_{\varepsilon,K}+
\mu^{\beta}_{\varepsilon,K}$
in place of $u$ and $v$, respectively, that
there exist  a  constant    $N \in(0,\infty)$,
depending only on $\tau$, the diameter of $\Omega$,   
    $d$, $K_{0}$,  $\bar H$, and $\delta$,
 and a  constant  $\eta>0$,
depending only on $\tau$,
    $d$,   and $\delta$,
such that, if $K\geq N\big(M^{\beta}_{\varepsilon,K}\big)^{ \eta}$, then \vspace{5pt} 
$$
u_{K}^{(\varepsilon)}-u_{-K}^{(\varepsilon)}
\leq  \mu^{\beta}_{\varepsilon,K}+2w^{\beta}_{\varepsilon,K}+
 NK^{-(\kappa-1)/4}+  NM^{\beta}_{\varepsilon,K}\omega\Big(
\big(M^{\beta}_{\varepsilon,K}\big)^{-1/\tau}K^{-1}\Big)
\vspace{5pt}$$
 in    $\Pi_{n_{0}-1}$, where 
$M^{\beta}_{\varepsilon,K}\geq1$ is any number satisfying  \vspace{5pt}
$$
M^{\beta}_{\varepsilon,K}\geq\big\|u_{K}^{(\varepsilon)}
-w^{\beta}_{\varepsilon,K},
u_{-K}^{(\varepsilon)}+w^{ \beta}_{\varepsilon,K}
+\mu^{\beta}_{\varepsilon,K}
\big\|_{C^{ \kappa }( \Pi_{n_{0}-1}  )}.
\vspace{5pt}$$

By letting $\beta\downarrow0$ we obviously obtain that,
if $K\geq N \big(M _{\varepsilon,K}\big) ^{ \eta}$, then  \vspace{5pt}
$$
u_{K}^{(\varepsilon)}-u_{-K}^{(\varepsilon)}
\leq  \mu _{\varepsilon,K}+2w _{\varepsilon,K}+     
 NK^{-(\kappa-1)/4}+  NM _{\varepsilon,K}\omega\big(
 M _{\varepsilon,K} ^{-1/\tau}K^{-1}\big)
\vspace{5pt}$$
 in    $\Pi_{n_{0}-1}$, where 
$$
\mu _{\varepsilon,K}=\sup_{\partial'
\Pi_{n_{0}-1}}\big(
u_{K}^{(\varepsilon)}-u_{-K}^{(\varepsilon)}-2w _{\varepsilon,K}\big)_{+},
\vspace{5pt}$$
and $M _{\varepsilon,K}\geq1$ is any number satisfying
$$
M _{\varepsilon,K}\geq\big\|u_{K}^{(\varepsilon)}
-w _{\varepsilon,K},
u_{-K}^{(\varepsilon)}+w _{\varepsilon,K}
+\mu _{\varepsilon,K}
\big\|_{C^{ \kappa }( \Pi_{n_{0}-1} )}.
\vspace{5pt}$$

First we discuss
what is happening as $\varepsilon\downarrow0$.  
By   
Theorem 1.1 of \cite{DKL_12}  
we obtain 
$w_{\varepsilon,K}\to 0$ in  
  $W^{1,2}_{p}(\Pi_{n_{0} })$
for any $p>1$,  
 which by
embedding theorems implies that $w_{\varepsilon,K}\to 0$ in
$C^{ \kappa }(\Pi_{n_{0} })$. Obviously, the constants 
$\mu_{\varepsilon,K}$ converge in 
$C^{ \kappa }(\Pi_{n_{0}-1})$ to  \vspace{5pt}   
$$ 
\sup_{\partial'
\Pi_{n_{0}-1}} (
u_{K} -u_{-K}   )_{+}.
\vspace{5pt}$$
Now 
Lemma \ref{lemma 2,4,1} (iv), applied in $\Pi_{n_{0}-1}$, implies that
for sufficiently small $\varepsilon$ one can take
$N\varepsilon ^{- \kappa }(n_{0})$ as $M_{\varepsilon,K}$,
where $\varepsilon(n_{0})$ is the distance
between the boundaries of $\Omega_{n_{0}}$
and $\Omega_{n_{0}-1}$ and $N$ is independent of $K$
 and $\varepsilon(n_{0})$. 
Thus, for sufficiently small $\varepsilon$,
if $K\geq N\varepsilon ^{- \kappa \eta}(n_{0})$,
then  \vspace{5pt}
$$
u_{K}^{(\varepsilon)}-u_{-K}^{(\varepsilon)}
\leq  \mu_{\varepsilon,K}+2w_{\varepsilon,K}+
 NK^{-(\kappa-1)/4}+N\varepsilon ^{- \kappa }(n_{0})
\omega\big(N\varepsilon ^{   \kappa /\tau}(n_{0})K^{-1}\big)
\vspace{5pt}$$
in $\Pi_{n_{0}-1}$, which after letting $\varepsilon\downarrow0$
yields  \vspace{5pt}
$$
u_{K} -u_{-K} 
\leq   
 NK^{-(\kappa-1)/4}+N\varepsilon ^{- \kappa }(n_{0})
\omega\big(N\varepsilon ^{  \kappa /\tau}(n_{0})K^{-1}\big)
+\sup_{\partial'
\Pi_{ n_{0}-1}} (
u_{K} -u_{-K}   )_{+}
\vspace{5pt}$$
in $\Pi_{n_{0}-1}$. Hence in $\Pi$
$$
u_{K}-u_{-K}\leq  NK^{-(\kappa-1)/4}
+N\varepsilon ^{- \kappa }(n_{0})
\omega\big(N\varepsilon ^{  \kappa /\tau}(n_{0})K^{-1}\big)
+\sup_{ \Pi\setminus 
\Pi_{n_{0}-1}} (
u_{K} -u_{-K}   )_{+}
\vspace{5pt}$$
$$
\leq  NK^{-(\kappa-1)/4}
+N\varepsilon ^{- \kappa }(n_{0})
\omega\big(N\varepsilon ^{  \kappa /\tau}(n_{0})K^{-1}\big)
+\xi(n_{0})
\medskip$$
where $\xi(n_{0})\to 0$ as $ n_{0}\to\infty$
by Lemma \ref{lemma 2,4,1} (iii).
This obviously proves the 
lemma  because, as is noted in Lemma \ref{lemma 2,4,1}
(ii),
  we have $u_{-K}\leq u_{K}$. \qed

\mysection{Existence 
of maximal and minimal $L_{p}$-viscosity
solutions}
                                                \label{section 12,2.1}

Fix   constants, $K_{0},  T \in(0,\infty)$, $p>d+2$,   
$\delta\in(0,1]$, and  fix a nonnegative
$$
\bar H\in L_{p}(\bR^{d+1}).
$$
Also according to the setting in Section 
 \ref{section 10.25.1}     
 we take a bounded domain $\Omega\subset\bR^{d}$
of class $C^{1,1}$ and set $\Pi=[0,T)\times\Omega$.

\begin{assumption}
                                    \label{assumption 12.2.1}
(i) The function $H$ is a nonincreasing 
function of $\sfu'_{0}$, is continuous with respect to 
  $\sfu'_{0}$,  uniformly
with respect to other variables
 $[\sfu'], (t,x)\in\bR^{d+1}$,
$\sfu''\in\bS$,
is measurable with respect to $(t,x)$ for any $\sfu$, and is
  Lipschitz continuous in $[\sfu']$ with Lipschitz constant
independent of $\sfu'_{0},\sfu'',(t,x)$.

(ii) For any $\sfu',(t,x)\in\bR^{d+1}$, the function $H(\sfu,t,x)$
is Lipschitz continuous with respect to $\sfu''$ and
at all points of differentiability of 
$H(\sfu,t,x)$  with
respect to $\sfu''$, we have
 $
D_{\sfu''}H  \in \bS_{\delta}$.

(iii) for all $\sfu',(t,x)\in\bR^{d+1}$
$$
 \big|H (\sfu',0,t,x)\big|\leq K_{0}|\sfu'|+\bar H(t,x).
\vspace{5pt}$$
 
\end{assumption}

\begin{assumption}
                                 \label{assumption 1,3.4}
We are given a function
$g\in C(\overline{\partial'\Pi})$.

\end{assumption}

We are going to use  
the following local version of Theorem 1.14 of \cite{Kr_17.1},
proved there for $g\in W^{1,2}_{p}(\bR^{d+1})$ with the solution
in global rather than local spaces $W^{1,2}_{p,\loc}$.
This local version is easier to prove because no boundary estimates
are needed and we will provide the proof elsewhere.

\begin{theorem}
                                       \label{theorem 10.9.1}
There exists a convex positive homogeneous of degree
one function $P(\sfu'')$ such that 
   at all points of its differentiability 
$D_{\sfu''}P \in
\bS_{\bar\delta}$, where $\bar\delta=\bar\delta(d,\delta)\in(0,\delta)$,
and for $P[u]=P(D^{2}u)$
and any $K>0$ there exists
  $v\in W^{1,2}_{p }(\Pi^{\rho})\cap C(\bar\Pi)$, 
 for any $\rho>0$,
such that $v=g$ on $\partial'\Pi$ and the equation  \vspace{5pt}     
\begin{equation}
                                              \label{1,31,4}
\partial_{t}v+\max\big(H[v],P [v]-K\big)=0,
\vspace{5pt}\end{equation}
holds \(a.e.\) in $\Pi$.

\end{theorem}

  By the maximum principle 
 the solutions $v=v_{K}$  
 are unique
and decrease as $K\to\infty$.

\begin{theorem}
                                          \label{theorem 1,31,3}
 Under the above assumptions, 
as $K\to\infty$, $v_{K}$ converges uniformly on
$\bar\Pi$ to a continuous function $v$ which
is an $L_{d+1}$-viscosity solutions of 
\eqref{7.29.10}
with  boundary
condition $v=g$ on $\partial'\Pi$.
Furthermore,  $v$ is 
the maximal $L_{d+1}$-viscosity subsolution of \eqref{7.29.10}
of class $C(\bar\Pi)$ with given boundary condition.
\end{theorem}

\begin{remark}
                                                  \label{remark 2,19,7}
To obtain an $L_{d+1}$-viscosity solution which is a minimal 
$L_{d+1}$-viscosity supersolution, 
it suffices to consider   \vspace{5pt}
\begin{equation}
                                             \label{3,18,1}
\partial_{t}v+\min\big(H[v],-P [-v]+K\big)=0,
\vspace{5pt}\end{equation}
which reduces to \eqref{1,31,4} if we replace $v$ with $-v$
and $H(\sfu,t,x)$ with $-H(-\sfu,t,x)$.

\end{remark}

This yields the following result.

\begin{theorem}
                                             \label{theorem 3,11,1}
Let $v_{-K}\in W^{1,2}_{p,\loc}(\Pi)\cap C(\bar\Pi)$
denote a unique solution of  \eqref{3,18,1}
\(a.e.\) in $\Pi$ with boundary data $v_{-K}=g$ on $\partial'\Pi$.
Then, as $K\to\infty$, $v_{-K}$ converges uniformly on
$\bar\Pi$ to a continuous function $w$ which
is an $L_{d+1}$-viscosity solutions of 
\eqref{7.29.10}
with  boundary
condition $w=g$ on $\partial'\Pi$.
Furthermore,  $w$ is 
the minimal $L_{d+1}$-viscosity supersolution of \eqref{7.29.10}
of class $C(\bar\Pi)$ with given boundary condition.

\end{theorem}

\begin{remark}
                                                  \label{remark 2,19,6} 
The existence of extremal $C$-viscosity solutions
is proved in the elliptic and parabolic cases 
in \cite{CKLS_99} when $H$ is a continuous function.
Our function $H(\sfu,t,x)$ is just measurable in $(t,x)$
and we are dealing with $L_{d+1}$-viscosity solutions.

Also note that the existence of
the extremal $L_{p}$-viscosity
solution for the elliptic case was proved in \cite{JS05}
with no continuity assumption on $H$ with respect to $x$.
We provide a method which in principle allows one to find it.

\end{remark}

Here is a stability result for the
extremal $L_{d+1}$-viscosity solutions.
In the following assumption there are two
 objects: $\kappa_{1}=\kappa(d,\delta)\in(1,2)$ (close to $1$),
  and 
 $\theta=\theta(\kappa,d,\delta)\in(0,1]$ (close to $0$),
$\kappa\in(1,\kappa_{1})$.
The values of $\kappa_{1}$ and $\theta$ are specified  
in the proof of Lemma 5.3 of \cite{Kr_14.1}.

 \begin{assumption}
                                  \label{assumption 10.23.1}
We have a representation 
$$
H(u,t,x)=F(u'',t,x)+G(u,t,x).
$$ 

(i) The functions $F$ and $G$ are measurable functions 
of their arguments.

(ii) For all values of the arguments
$$
|G(u,t,x)|\leq K_{0}|u'|+\bar{H}(t,x).
$$
 
(iii) The function $F$ is positive homogeneous
of degree one with respect to $u''$,
is Lipschitz continuous with respect  to
$u''$, and at all points of differentiability of
$F$ with respect to $u''$ we have $D_{\sfu''}F \in\bS_{\delta}$.

(iv) For any $R\in(0,R_{0}]$, $(t,x)\in\bR^{d+1}$,
 and $u''\in\bS$
with $|u''|=1$ ($|u''|:=(\tr u''u'')^{1/2}$), we have
$$
\theta_{R,t,x}:=\dashint_{C_{R}(t,x)} 
 |F(u'',s,y)-\bar{F}_{R, x}(u'',s)|\,dsdy\leq\theta,
$$
where
$$
\bar{F}_{R, x}(u'',s)=\dashint_{B_{R}(x)}F(u'',s,y)\,dy.
$$

\end{assumption}
Take 
$\kappa_{1}= \kappa_{1}(d,\delta)\in(1,2)$ and
$ \theta=\theta(\kappa,d,\delta)\in(0,1]$ that are
introduced before  Assumption 
 \ref{assumption 10.23.1}. Then fix a $\kappa$
satisfying
$$
1<\kappa<[2-(d+2)/p]\wedge \kappa_{1}.
$$

\begin{theorem}
                                             \label{theorem 8,12,1}
Let $H_{n}$, $n=0,1,...$,  satisfy Assumptions
\ref{assumption 12.2.1} and \ref{assumption 10.23.1}   
with   
$R_{0}\in(0,1]$ independent of $n$, with
 $K_{0}$ and $\bar H$ from the beginning of the section,
     $\theta(\kappa,d,\bar\delta)/2$
in place of $\theta$ \(specified above\) and have, perhaps,
different Lipschitz constants with respect to $[\sfu']$
for different $n$.
Suppose that $H_{0}$ is Lipschitz continuous in $\sfu$
with a constant independent of $(t,x)$.
Let $v_{n}$, $n=0,1,...$, 
be the maximal $L_{d+1}$-viscosity  solutions  
  of class $C(\bar\Pi)$ of
$\partial_{t}v_{n}+H_{n}[v_{n}]=0$  
in $\Pi$ with boundary condition $v_{n}=g_{n}$ on $\partial'\Pi$,
where $g_{n}\in C(\bar\Pi)$ and $g_{n}\to g_{0}$ in
$C(\bar\Pi)$ as $n\to\infty$.

Assume that for any $M>0$
\begin{equation}
                                             \label{8,12,7}
\Delta_{n,M}(t,x):=\sup_{|\sfu|\leq M}
\big|H_{n}(\sfu,t,x)-H_{0}(\sfu,t,x)\big|\to0
\end{equation}
in $L_{d+1}(\Pi)$ as $n\to\infty$. Also
assume that for all values of the arguments and $n$
\begin{equation}
                                                   \label{8,12,1}
\big|H_{n}(\sfu,t,x)-H_{0}(\sfu,t,x)\big|\leq
\bar H(t,x)\big(1+|\sfu'|\big).
\end{equation}

Then $v_{n}\to v_{0}$ in
$C(\bar\Pi)$ as $n\to\infty$.
The same holds true if $v_{n}$ are minimal
$L_{d+1}$-viscosity solutions of class $C(\bar\Pi)$.

\end{theorem}

Proof.  According to Theorem \ref{theorem 1,31,3},   it suffices to show
that   \vspace{5pt}
\begin{equation}
                                                        \label{8,12,2}
\sup_{K\geq 1}\big(v_{n,K}-v_{0,K}\big)\to 0,
\vspace{5pt}\end{equation}
in $C(\bar\Pi)$, where $v_{n,K}$ are the solutions of    \vspace{5pt}
$$
\partial_{t}v_{n,K}+\max\big(H_{n}[v_{n,K}], P[v_{n,K}]-K\big)=0
\vspace{5pt}$$
in $\Pi$ (a.e.) 
with boundary condition $v_{n,K}=g_{n}$ on $\partial'\Pi$. 

Observe that
$$
\big|\max\big(H_{n}[v_{n,K}], P[v_{n,K}]-K\big)-
\max\big(F_{n}[v_{n,K}], P[v_{n,K}]-K\big)\big|
$$
$$
\leq |G_{n}[v_{n,K}]|=|b^{i} D_{i}v_{n,K}+cv_{n,K}+\tau \bar H|
$$
for certain    
measurable $\bR^{d}$-valued $ b $, and measurable
real-valued $\tau $ and $c$, such that
$| b |, |c|\leq K_{0}$, and $|\tau|\leq1$.
Also $\max\big(F_{n}[0], P[0]-K\big)=0$.
Therefore, 
by the mean-value theorem we have    \smallskip   
\begin{equation}
                                           \label{1,9,81}
\partial_{t}v_{n,K}+
a^{ij}D_{ij}v_{n,K}  +b^{i} D_{i}v_{n,K}+cv_{n,K}+\tau \bar H=0
\end{equation}

\noindent
(a.e.) in $\Pi$ for certain  measurable 
$\bS_{ \bar\delta}$-valued $(a^{ij} )$, and perhaps different
measurable $\bR^{d}$-valued $ b $, and measurable
real-valued $\tau $ and $c$, such that
$| b |, |c|\leq K_{0}$, and $|\tau|\leq1$.
By the parabolic Aleksandrov estimates
\eqref{1,9,81} implies that $|v_{n,K}|$
are uniformly bounded in $\bar\Pi$ and
  by the linear theory of parabolic equations we conclude that
the family $\{v_{n,K}:K\geq1, n\geq0\}$ is precompact in 
$C(\bar\Pi)$.

Next, fix a $\rho>0$ such that $\Pi^{\rho}\ne\emptyset$
and observe  that,  as we know from
\cite{DKL_12}, \cite{Kr_12.3},
there is a number $\gamma=\gamma(d,\delta,K_{0})\in(0,1)$ such that
 there is a constant $N$,
depending only on $\rho$, $d$, $\delta$, and $K_{0}$,
such that for any cylinder $C_{\rho}(t_{0},x_{0})\subset\Pi$
  we have due to \eqref{1,9,81} that \smallskip
$$
\int_{C_{\rho}(t_{0},x_{0})}\big|D^{2}v_{n,K}\big|^{\gamma}
\,dxdt\leq N\sup_{\Pi}|v_{n,K}|+
N\Big(\int_{\Pi}|cv_{n,K}+\tau \bar H|^{d+1}\,dxdt\Big)^{\gamma/(d+1)}
\medskip$$
for all $n,K$. Here the right-hand side is dominated by a constant
independent of $n,K$, and it follows, by Chebyshev's inequality
 that   there is 
a constant $N$ (perhaps depending on $\rho$)
for which  \smallskip
\begin{equation}
                                                          \label{8,12,3}
\big|\Pi^{\rho}\cap\big\{|D^{2}v_{n,K}|\geq M\big\}\big|
\leq N  M^{-\gamma} 
\medskip\end{equation}
 for all $n\geq0,K\geq1, M>0$.

To finish with preparations, 
set 
$$
H_{n,K}=\max(H_{n},P-K),\quad
 F_{n,K}=\max (F_{n},P),
$$
$ G_{n,K}=H_{n,K}-F_{n,K} $, where $F_{n}$ is taken from
Assumption \ref{assumption 10.23.1} written for $H_{n}$.
Then $F_{n,K}$ and $G_{n,K}$ satisfy Assumption \ref{assumption 10.23.1}
(i) and (ii) with the same $K_{0}$ and $\bar H$.
Assumption \ref{assumption 10.23.1} (iii) also is satisfied with     
  $\bar\delta$ in place of $\delta$. Finally, easy
manipulations, using the fact that in the assumptions
of the theorem we suppose that Assumption \ref{assumption 10.23.1}
(iv) is satisfied for $F_{n}$
with $\theta(\kappa,d, \bar\delta )/2$ in place of $\theta$,
show that Assumption \ref{assumption 10.23.1}
(iv) is satisfied for $F_{n,K}$ with 
$\theta=\theta(\kappa,d, \bar\delta) $.

Thus, all the assumptions of Theorem 2.1 of \cite{Kr_14.1}  
are satisfied apart from $g\in W^{1,2}_{\infty}
(\bR^{d+1})$  and $\Omega\in C^{2}$. We will show in a 
separate publication that these assumptions can be replaced
with the current ones.
Now
  since $v_{n,K}$ is a classical solution   of \eqref{1,31,4},
we obtain from that theorem, for any small $\rho>0$,
the estimates of the 
$C^{1+\alpha}(\Pi^{\rho})$-norms of $v_{n,K}$ uniform
with respect to $n$ and $K$.
  Therefore,
by   interpolation   
theorems  we get \vspace{5pt}
\begin{equation}
                                                          \label{8,12,4}
\sup_{\Pi^{\rho}}|Dv_{n,K}|\leq N ,
\vspace{5pt}\end{equation}
where and below by $N$ we denote various constants independent
of $K$ and $n$, perhaps depending on $\rho$.

Now set $w_{n,K}=v_{0,K}-v_{n,K}$, and observe that   \vspace{5pt}
$$
0=\partial_{t}w_{n,K}+I_{1}+I_{2}+I_{3},
$$
where
$$
I_{1}=\max\big(H_{0}[v_{0,K}], P[v_{0,K}]-K\big)
\vspace{5pt}$$
$$
-
\max\Big(H_{0}\big(v_{0,K},Dv_{0,K},D^{2}v_{n,K}
\big), P[v_{n,K}]-K\Big)=a^{ij}D_{ij}w_{n},
\vspace{5pt}$$
$$
I_{2}=\max\Big(H_{0}\big(v_{0,K},Dv_{0,K},D^{2}v_{n,K}\big), 
P[v_{n,K}]-K\Big)
\vspace{5pt}$$
$$
-\max\Big(H_{0}\big(v_{n,K},Dv_{n,K},D^{2}v_{n,K}\big), 
P[v_{n,K}]-K\Big),
\vspace{5pt}$$
$$
I_{3}=\max\big(H_{0}[v_{n,K}], P[v_{n,K}]-K\big)
\vspace{5pt}$$
$$
-\max\big(H_{n}[v_{n,K}], P[v_{n,K}]-K\big),
\vspace{5pt}$$
and $(a^{ij})$ is an $\bS_{ \bar\delta}$-valued function. By assumption \vspace{5pt}
$$
 |I_{2}| \leq  N\big(
|w_{n,K}|+|Dw_{n,K}|\big).
\vspace{5pt}$$
Therefore,
$$
0=\partial_{t}w_{n,K}+a^{ij}D_{ij}w_{n,K}+b^{i}D_{i}w_{n,K}+
cw_{n,K}+I_{3},
$$
where $b$ and $c$ are bounded uniformly with respect to $n,K$.

Upon observing that
 by assumption and \eqref{8,12,4} in $\Pi^{\rho}$ for any $M>0$ we have
$$
|I_{3}| \leq \bar HN I_{|D^{2}v_{n,K}|\geq M}
+\Delta_{n,M+N }I_{|D^{2}v_{n,K}|\leq M} 
$$
and using the parabolic Aleksandrov estimates in $\Pi^{\rho}$
 we conclude that
  there exists a constant $N $ such
that for all $n,K,M$ in $\Pi$
\begin{align}
|v_{0,K}-v_{n,K}|\leq N \|\bar H I_{|D^{2}v_{n,K}|\geq M}
\|_{L_{d+1}(\Pi^{\rho})}
\nonumber\\[10pt]
                                              \label{8,12,6}
+N \|\Delta_{n,M+N }
\|_{L_{d+1}(\Pi^{\rho})}
+\sup_{\Pi\setminus \Pi^{\rho}} |v_{0,K}-v_{n,K}|.
\end{align}
Here
$$
\sup_{n,K}\|\bar H I_{|D^{2}v_{n,K}|\geq M}\|_{L_{d+1}(\Pi^{\rho})}
\to 0 
$$
as $M\to \infty$, since $\bar H\in L_{d+1}(\Pi)$
and \eqref{8,12,3} holds. Therefore, by first taking
the sup's with respect to $K\geq1$ in \eqref{8,12,6},
then sending $n\to\infty$, using assumption \eqref{8,12,7},
and then sending $M\to\infty$, we infer from
\eqref{8,12,6} that, for any small $\rho>0$
 $$
\nlimsup_{n\to\infty}\sup_{\Pi}|v_{0,K}-v_{n,K}|
\leq
\nlimsup_{n\to\infty}\sup_{\Pi\setminus\Pi^{\rho}}|v_{0,K}-v_{n,K}|.
$$
After that it only remains to set $\rho\downarrow0$
and   use the equicontinuity
of $v_{n,K}$ and the fact that $g_{n}\to g_{0}$
uniformly in $\overline{\partial'\Pi}$.
The theorem is proved.  \qed

\begin{remark}
                                       \label{remark 8,14,1}
It follows from the above proof that $\bar H$ in
\eqref{8,12,1} can be replaced with $\bar H_{n}$,
provided that the family $|\bar H_{n}|^{d+1}$
is uniformly integrable over $\Pi$.
\end{remark}

An obvious consequence of this theorem
is the stability of uniqueness.

\begin{corollary}
                                               \label{corollary 8,12,1}
Suppose that for any $n=1,2,...$ there is only one
  $L_{d+1}$-viscosity  solutions  
  of class $C(\bar\Pi)$ of
$\partial_{t}v_{n}+H_{n}[v_{n}]=0$  
in $\Pi$ with boundary condition $v_{n}=g_{n}$ on $\partial'\Pi$.
Then the same holds for $n=0$.
\end{corollary}

Coming back to Theorem \ref{theorem 1,31,3}, observe that,
as we have mentioned above,
by the maximum principle $v_{K}$ decreases as $K$ increases.  
The precompactness of $\{v_{K},K\geq1\}$ in $C(\bar\Pi)$ is
proved in the same way as in the above proof after
\eqref{1,9,81} using the fact that
$$
 \max\big(H_{n}[v_{ K}], P[v_{ K}]-K\big)-
\max\big(H_{n}(v_{ K},Dv_{K},0,\cdot), P[0]-K\big)=a^{ij}D_{ij}v_{K},
$$
where $(a^{ij})$ is $\bS_{\bar\delta}$-valued and
$$
\max\big(H_{n}(v_{ K},Dv_{K},0,\cdot), P[0]-K\big)=
\max\big(H_{n}(v_{ K},Dv_{K},0,\cdot),-K\big)
$$
$$
= b^{i} D_{i}v_{n,K}+cv_{n,K}+\tau \bar H.
$$

It follows   that $v_{K}$ converges uniformly on $\bar\Pi$
 as $K\to\infty$
to a function $v\in C(\bar\Pi)$. To prove that
$v$ is an $L_{d+1}$-viscosity solution we need the following.

\begin{lemma}
                                           \label{lemma 2,2,2}
There is a constant $N$ depending only on $d$,
$\delta$, and the Lipschitz constant of $H$ with respect to
$[\sfu']=\big(\sfu'_{1},...,\sfu'_{d}\big)$ such that
for any $r\in(0,1]$ and $C_{r}(t,x)$ satisfying $C_{r}(t,x)\subset\Pi$  and
$\phi\in W^{1,2}_{d+1}\big(C_{r}(t,x)\big)$ we have on $C_{r}(t,x)$ that
\begin{equation}
                                             \label{09.20.1}
v\leq \phi+Nr^{d/(d+1)}\big\|(\partial_{t}\phi+
H[\phi])^{+}\big\|_{L_{d+1}(C_{r}(t,x))}
+\max_{\partial'C_{r}(t,x)}(v-\phi)^{+} .
\vspace{5pt}\end{equation}
\begin{equation}
                                             \label{09.20.2}
v\geq \phi-Nr^{d/(d+1)}\big\|(\partial_{t}\phi+
H[\phi])^{-}\big\|_{L_{d+1}(C_{r}(t,x))}
-\max_{\partial'C_{r}(t,x)}(v-\phi)^{-} .
\vspace{5pt}\end{equation}
\end{lemma}

Proof. Observe that  in $C_{r}(t,x)$ (a.e.) \vspace{5pt}
$$
-\partial_{t}\phi-\max\big(H[\phi],P[\phi]-K\big)=
-\partial_{t}\phi-\max\big(H[\phi],P[\phi]-K\big)
\vspace{5pt}$$
$$
+\partial_{t}v_{K}+
\max\big(H[v_{K}],P[v_{K}]-K\big)
\vspace{5pt}$$
$$
=
\partial_{t} (v_{K}-\phi )+a_{ij}D_{ij}(v_{K}-\phi )
+b_{i}D_{i}(v_{K}-\phi )-c(v_{K}-\phi ),
\vspace{5pt}$$
where $a=(a_{ij})$ is an $\bS_{\bar\delta}$-valued
function,
$b=(b_{i})$ is bounded by the Lipschitz
constant of $H$ with respect to $[\sfu']$, and $c\geq0$.
It follows by  the parabolic
Aleksandrov estimates that  \vspace{5pt}
$$
v_{K}\leq \phi+\max_{\partial'C_{r}(t,x)}(v_{K}-\phi)^{+}
\vspace{5pt}$$
\begin{equation}
                                             \label{09.20.3}
+Nr^{d/(d+1)}\big\|\big(\partial_{t}\phi+
\max\big(H[\phi],P[\phi]-K\big)\big)^{+}\big\|_{L_{d+1}(C_{r}(t,x))},
\vspace{5pt}\end{equation}
where the constant $N$ is of the type described
in the statement of the present lemma.
 We obtain \eqref{09.20.1} from \eqref{09.20.3}
by letting $K\to\infty$.
In the same way \eqref{09.20.2} is established.
The lemma is proved.     \qed

{\bf Proof of Theorem \ref{theorem 1,31,3}}.   
First we  prove that $v$ is  an $L_{d+1}$-viscosity solution.   
Let $(t_{0},x_{0})\in\Pi$ and $\phi\in W^{1,2}_{d+1,\loc}
(\Pi)$ be such that $v-\phi$ attains a local maximum at $(t_{0},x_{0})$.
Then for   $\varepsilon
>0$ and all small $r>0$ for
$$
\phi_{\varepsilon,r}(t,x)=\phi (t,x)-\phi(t_{0},x_{0})
+v(t_{0},x_{0})+\varepsilon\big(
|x-x_{0}|^{2}+t-t_{0}- r^{2}\big)
$$
 we have that 
$$
\max_{\partial'C_{r}(t_{0},x_{0})}(v -\phi_{\varepsilon,r} )^{+}
=0.
\vspace{5pt}$$
Hence, by Lemma \ref{lemma 2,2,2}
\begin{align*}
   \varepsilon r^{2}= 
(v -\phi_{\varepsilon,r})(t_{0},x_{0})
\leq Nr^{d/(d+1)}\big\|\big(\partial_{t}\phi_{\varepsilon,r}+
H[\phi_{\varepsilon,r}]\big)^{+}\big\|_{L_{d+1}(C_{r}(t_{0},x_{0}))},
\\[10pt]
\varepsilon^{d+1}\leq Nr^{-(d+2)}\big\|\big(\partial_{t}\phi_{\varepsilon,r}+
H[\phi_{\varepsilon,r}]\big)^{+}\big\|^{d+1}_{L_{d+1}(C_{r}(t_{0},x_{0}))}.
\end{align*}
By letting $r\downarrow0$ and using the continuity
of $H(\sfu,t,x)$ in $\sfu'_{0}$, which is assumed to be uniform
with respect to other variables and also using the continuity of
$\phi$ (embedding theorems) and $v$, we obtain  \vspace{5pt}
\begin{equation}
                                        \label{6.8.40}
N\lim _{ r\downarrow0}\esssup_{C_{r}(t_{0},x_{0}) }
\Big(\partial_{t}\phi_{\varepsilon }(t,x)+
H\big(v(t,x),D\phi_{\varepsilon }(t,x),
D^{2}\phi_{\varepsilon }(t,x)\big)\Big)\geq \varepsilon.
\vspace{5pt}\end{equation}
where $\phi_{\varepsilon }=\phi+\varepsilon
\big(|x-x_{0}|^{2}+t-t_{0}\big)$.  Finally, observe that
  $H(\sfu,t,x)$ is Lipschitz continuous with respect to   
$\big([\sfu'],\sfu''\big)$  with Lipschitz constant
independent of $\sfu'_{0} ,(t,x)$ by assumption. 
Then letting $\varepsilon\downarrow0$ in \eqref{6.8.40}
proves that
$v$ is an $L_{d+1}$-viscosity subsolution.
The fact that it is also an $L_{d+1}$-viscosity supersolution
is proved similarly on the basis of \eqref{09.20.2}.

Finally, we prove that $v$ is the maximal continuous
$L_{d+1}$-viscosity subsolution. Let $w$ be an 
$L_{d+1}$-viscosity subsolution of \eqref{7.29.10}
 of class 
$C(\bar \Pi) $ with boundary data $g$. To prove that $v\geq w$,
it suffices to show that for any $\varepsilon>0$
and $K>1$ we have $u_{K}+\varepsilon(T-t)\geq w$
in $\bar\Pi$.

Assume the contrary and observe that, since
$u_{K}+\varepsilon(T-t)-w\geq0$ on $\partial'\Pi$,
there is a point $(t_{0},x_{0})\in\Pi$ such that  \vspace{5pt}
\begin{equation}
                                                 \label{3,19.1}
\gamma:= u_{K}(t_{0},x_{0})+\varepsilon(T-t_{0})-w(t_{0},x_{0})<0
\vspace{5pt}\end{equation}
and $u_{K}+\varepsilon(T-t)- w\geq \gamma$ in $C_{\rho}(t_{0},x_{0})$
for sufficiently small $\rho>0$,
so that by definition    \vspace{5pt}
$$
-\varepsilon+\lim _{ r\downarrow0}\esssup_{C_{r}(t_{0},x_{0}) }
\Big(\partial_{t}u_{K}(t,x)+
H\big(w(t,x),Du_{K}(t,x),
D^{2}u_{K}(t,x),t,x\big)\Big)\geq 0.
\vspace{5pt}$$
Since $H$ is a decreasing function of $\sfu'_{0}$, in light of
\eqref{3,19.1},    \vspace{5pt}
$$
-\varepsilon+\lim _{ r\downarrow0}\esssup_{C_{r}(t_{0},x_{0}) }
\Big(\partial_{t}u_{K}(t,x)+
H\big(u_{K}(t,x),Du_{K}(t,x),
D^{2}u_{K}(t,x),t,x\big)\Big)\geq 0.
\vspace{5pt}$$
This is however impossible since $\partial_{t}u_{K}
+H[u_{K}]\leq0$ in $\Pi$ (a.e.).
This contradiction finishes proving the theorem.   \qed

\mysection{Uniqueness 
of $L_{p}$-viscosity solutions
for parabolic Isaacs \\ \quad\, equations}
                                               \label{section 3,7,1}

Fix some constants $\delta\in(0,1]$, $K_{0},T\in(0,\infty)$, $p>d+2$.
Assume that we are given countable sets
  $A$ and $B$,   and,
for each $\alpha\in A$ and $\beta\in B$, we are given
an $\bS_{\delta} $-valued
   function $a^{\alpha\beta} $ on $\bR^{d+1}=\big\{(t,x):t\in\bR,x\in\bR^{d}\big\}$,
a real-valued function $b^{\alpha\beta}(\sfu',t,x)$
on $ \bR^{d+1}\times \bR^{d+1}$, and a real-valued $\bar H(t,x)\geq0$
on $\bR^{d+1}$.

\begin{assumption}
                                   \label{assumption 8,13,1}
\(i\,\) The above introduced functions are measurable.

\(ii\,\) The function $a^{\alpha\beta}(t, x )$ is uniformly continuous
with respect to $(t,x)$ uniformly with respect to $\alpha,\beta$ and,
with $\gamma$ introduced before Assumption
\ref{assumption 3,13,1}, for all values of indices and arguments
$$
\big|a^{\alpha\beta}(t,x)-a^{\alpha\beta}(t,y)\big|\leq K_{0}|x-y|^{\gamma}.
$$

\(iii\,\) The function $b^{\alpha\beta}(\sfu',t,x)$ is nonincreasing
with respect to $\sfu'_{0}$, is Lipschitz continuous
with respect to $\sfu'$ with Lipschitz constant $K_{0}$, and
for all values of  indices and  arguments
$$
\big|b^{\alpha\beta}(\sfu',t,x)\big|\leq  K_{0}|\sfu'| +\bar H(t,x).
$$

 (iv)  We have $\bar H\in L_{p}(\bR^{d+1})$.    

\end{assumption}
 
For
$$
\sfu=(\sfu',\sfu''),\quad \sfu'=\big(\sfu'_{0},\sfu'_{1},...,\sfu'_{d}\big)
\in \bR^{d+1},
\quad \sfu''\in\bS,
\vspace{5pt}$$
and $(t,x)\in \bR^{d+1}$
introduce   \vspace{5pt}
$$
H(\sfu,t,x)=
\supinf_{\alpha\in A\,\,\,\beta\in B}
\big[a^{\alpha\beta}_{ij}( t,x)\sfu''_{ij}  +
b ^{\alpha\beta} (\sfu',t, x) \big] 
\vspace{5pt}$$
(where as everywhere the summation convention is enforced
and the summations are done inside the brackets).

For  sufficiently smooth functions $u=u(t,x)$ define  \vspace{5pt} 
\begin{equation}
                                                      \label{1.16.1}
H[u](t,x)=  H\big(u(t,x),Du(t,x),D^{2}u(t,x) ,t,x\big).
\vspace{5pt}\end{equation}
 
As is usual in this article we take
 an open bounded subset $\Omega$ of $\bR^{d}$ of class $C^{1,1}$
 and  set
$$
\Pi=[0,T)\times\Omega.
$$

Here is   the main results of this section.
\begin{theorem}
                                                   \label{theorem 2,3,1}
Under the above assumption for any $g\in C(\overline{\partial'\Pi})$
 there exists
  a unique continuous in $\bar \Pi$,
$L_{d+1}$-viscosity solution $g$ of the Isaacs equation
\begin{equation}
                                                     \label{4.1.4}
\partial_{t}u+H[u]=0 
\vspace{5pt}\end{equation}
in $\Pi$  with boundary condition $u=g$ on
$\partial' \Pi$. 

\end{theorem}

\begin{remark}
                                     \label{remark 2,19,8}
Under Assumption \ref{assumption 8,13,1} (i), (ii),
in case 
$$
b^{\alpha\beta}(\sfu',t,x)=\sum_{i=1}^{d}b^{\alpha\beta}_{i}( t,x)\sfu'_{i}
-c^{\alpha\beta}( t,x)\sfu'_{0}+f^{\alpha\beta}( t,x)
$$
 with uniformly continuous (in $(t,x)$
uniformly in $\alpha,\beta$) and uniformly bounded
 coefficients and the free terms,
the uniqueness of $C$-viscosity   
solutions  
 is stated without proof
in  Theorem 9.3 in \cite{CKS00}. This case
is covered by Theorem \ref{theorem 3,13,1}.

For general $H$, not necessarily related
to Isaacs equations, uniqueness is claimed for $\partial_{t}u  
+H(D^{2}u)=0$ in Lemma 4.7 of \cite{Wa92_1}.
It is proved for $L_{p}$-viscosity solutions
 in Lemma 6.2 of  \cite{CKS00}
in case $H$ is independent of $(t,x)$ with no reference to Wang's 
Lemma 4.7 of \cite{Wa92_1}.

In the elliptic case Jensen and  \'Swi{\c e}ch  \cite{JS05}  
proved the uniqueness
of continuous $L_{p}$-viscosity solutions for Isaacs equations,
assuming that $b^{\alpha\beta}_{i}$, $c^{\alpha\beta}$
are bounded,   $\sup_{\alpha,\beta}|f^{\alpha\beta}|
\in L_{p}$ and   Assumptions \ref{assumption 8,13,1} (i), (ii)  
are satisfied. Their proof uses
a remarkable Corollary 1.6 of
 \'Swi{\c e}ch \cite{Sw_97}
of which the parabolic counterpart is given in \cite{CKS00}.

An important difference with \cite{JS05} here
is that we consider lower-order terms in a more general form,
but in \cite{JS05} the summability assumption is weaker
(some $p<d$ are allowed).

\end{remark}

We prove Theorem \ref{theorem 2,3,1}
by using Corollary \ref{corollary 8,12,1} after some
preparations.
 
Fix a nonnegative 
$\zeta\in C^{\infty}_{0}(\bR^{d+1})$ which integrates to one
and  for $n=1,2,...$,  introduce $\zeta_{n}(t,x)
=n^{ d+1}\zeta (nt ,nx )$.  
Also, for real-valued $\xi=\xi^{\alpha\beta}$ given on 
$A\times B$ and $\sfu',(t,x)\in\bR^{d+1}$ define \smallskip
$$
\cH_{ 0 }(\xi,\sfu',t,x)=\supinf_{\alpha\in A\,\,\,\beta\in B}
\big[\xi^{\alpha\beta}  +
b ^{\alpha\beta} (\sfu',t, x) \big],
\vspace{5pt}$$
$$
\cH_{ n }(\xi,\sfu',t,x)=\cH_{ 0 }(\xi,\sfu',t,x)*\zeta_{n}(t,x),
\quad n=1,2,...,
\medskip$$
where the convolution is performed with respect to $(t,x)$.
We will also use the notation
$$
|\xi|=\supsup_{\alpha\in A\,\,\,\beta\in B}|\xi^{\alpha\beta}|.
$$

Here are a few properties of  $\cH^{ n }$.

\begin{lemma}
                                              \label{lemma 8,13,1}
\(i\,\) For each $n\geq1$, the function
$\cH_{ n }(\xi,\sfu',t,x)$ is continuous, is Lipschitz
continuous in $\sfu'$ with Lipschitz constant $K_{0}$,
and the function $\big|\cH_{ n }(0,0,t,x)\big|$ is bounded.

\(ii\,\) For each $n\geq1$, the function
$\cH_{ n }(\xi,\sfu',t,x)$ is infinitely differentiable
with respect to $x$ \(and $t$\), and there exists a constant $N$
\(depending on $n$\) such that $|D_{x}\cH_{ n }|\leq N \big(1+|\sfu'|\big)$
for all values of arguments.

\(iii\,\) For all values of the arguments
$$
\big|\cH_{ n }(\xi,\sfu',t,x)-\cH_{0}(\xi,\sfu',t,x)\big|
\leq 2K_{0}|\sfu'|+\bar H+\hat H,
$$
where
$$
\hat H=\sup_{n}\hat H*\zeta_{n}.
$$

\(iv\,\) For each $M>0$
\begin{equation}
                                                      \label{8,13,1}
\delta_{n,M}(t,x):=\sup_{|\sfu'|\leq M}
\sup_{|\xi|\leq M}\big|
\cH_{ n }(\xi,\sfu',t,x)-\cH_{0}(\xi,\sfu',t,x)\big|\to0
\end{equation}
in $L_{d+1}(\Pi)$ as $n\to\infty$.
\end{lemma}

Proof. Assertions (i)-(iii) are quite elementary
and their proofs are left to the reader. To prove
(iv) set
$$
\delta_{n}(\xi,\sfu',t,x)
=\cH_{ n }(\xi,\sfu',t,x)-\cH_{0}(\xi,\sfu',t,x)
$$
and observe that  for any $(\xi,\sfu')$
$$
\lim_{n\to\infty}\big\|
\delta_{n}(\xi,\sfu',\cdot) \big\|
_{L_{d+1}(\Pi)}\to0
$$
by the $L_{d+1}$-continuity of $L_{d+1}$-functions.
Furthermore, by the Lipschitz continuity
of $\cH_{n}$ with respect to $(\xi,\sfu')$
uniform with respect to $n,t,x$, for any
$\varepsilon>0$, one can find $m$
and $(\xi_{k},\sfu'_{k})$, $k=1,...,m$, such that
$|\xi_{k}|,|\sfu'_{k}|\leq M$ and any
$(\xi,\sfu')$ satisfying $|\xi |,|\sfu' |\leq M$
has a neighbor $(\xi_{k},\sfu'_{k})$ such that
$$
 \big|
\cH_{ n }(\xi,\sfu',t,x)-\cH_{n}(\xi_{k},\sfu'_{k},t,x)\big |
 \leq\varepsilon
$$
in $\Pi$ for any $n\geq0$. It follows that
$$
\delta_{n,M}(t,x)\leq \max_{k=1,...,m}
\big|
\delta_{n}(\xi_{k},\sfu'_{k},t,x) \big|
+\varepsilon.
$$
Then
$$
\nlimsup_{n\to\infty}\|\delta_{n,M}\|_{L_{d+1}(\Pi)}
\leq\nlimsup_{n\to\infty}\sum_{k=1}^{m}
\big\|
\delta_{n}(\xi_{k},\sfu'_{k},\cdot) \big\|
_{L_{d+1}(\Pi)}+ N\varepsilon= N\varepsilon,
$$
 where $N$ is independent of $\varepsilon$. This certainly proves (iv)
and the lemma.   \qed

{\bf Proof of Theorem \ref{theorem 2,3,1}}.
First we check that the assumptions of Theorem \ref{theorem 8,12,1}
are satisfied for
$$
H_{n}( \sfu ,t,x)=\cH_{n} (\xi,\sfu',t,x ),
\quad \xi^{\alpha\beta}=a^{\alpha\beta}_{ij}\sfu''_{ij}.
$$
Assumptions \eqref{8,12,7}, \eqref{8,12,1}
are taken care of in Lemma \ref{lemma 8,13,1}. Assumption
\ref{assumption 12.2.1}
 is easily seen to be satisfied with $\bar H+\hat H$
in place of $\bar H$.

In what concerns Assumption \ref{assumption 10.23.1}
we set
$$
F_{n}(\sfu',t,x)=F_{0}(\sfu',t,x)=\supinf_{\alpha\in A\,\,\,\beta\in B}
 a^{\alpha\beta}_{ij}( t,x)\sfu''_{ij} .
$$
Then Assumptions \ref{assumption 10.23.1} (ii), (iii)
are obviously satisfied. 
The remaining assumptions of Theorem \ref{theorem 8,12,1}
is Assumption \ref{assumption 10.23.1} (iv),
which is satisfied for any $\theta>0$ if $R_{0}$ is chosen
appropriately in light of Assumption \ref{assumption 8,13,1} (ii).

Thus Theorem \ref{theorem 8,12,1} is applicable.
Furthermore, for each $n\geq1$ the assumptions
of Theorem \ref{theorem 3,13,1} are satisfied
in light of Lemma \ref{lemma 8,13,1} (ii)
and Assumption \ref{assumption 8,13,1} (ii).
Also uniqueness of $C$-viscosity solutions
implies that of $L_{p}$-viscosity solutions.

Now, the combination of Theorem \ref{theorem 3,13,1}
and Corollary \ref{corollary 8,12,1}
(proved under the assumptions of Theorem
\ref{theorem 8,12,1}) immediately yields the desired result. \qed

\end{document}